\documentclass[11pt,a4paper]{article}

\usepackage{epsf,epsfig,amsfonts,amsgen,amsmath,amstext,amsbsy,amsopn,amsthm
}
\usepackage{ebezier,eepic}
\usepackage{color}
\usepackage{multirow}
\newtheorem{theorem}{Theorem}

\newtheorem{lemma}[theorem]{Lemma}
\newtheorem{false statement}{False statement}

\theoremstyle{definition}

\newtheorem{claim}{Claim}

\newenvironment{pf}{\noindent{\em Proof.}}
{\enspace\vrule height5pt depth0pt width5pt}

\newcounter{mathitem}
  {\begin{list}{{$(\roman{mathitem})$}}{
   \setcounter{mathitem}{0}
   \usecounter{mathitem}
   \setlength{\topsep}{0pt plus 2pt minus 0pt}
   \setlength{\parskip}{0pt plus 2pt minus 0pt}
   \setlength{\partopsep}{0pt plus 2pt minus 0pt}
   \setlength{\parsep}{0pt plus 2pt minus 0pt}
   \setlength{\leftmargin}{35pt}
   \setlength{\itemsep}{0pt plus 2pt minus 0pt}}}
  {\end{list}}

\begin{document}

\title{\bf\Large Coloring graphs with two odd cycle lengths}

\date{}

\author{Jie Ma\footnote{School of Mathematical Sciences,
University of Science and Technology of China, Hefei, 230026,
P.R. China. Email: jiema@ustc.edu.cn. Partially supported by NSFC grants 11501539 and 11622110.}~
~~~~~~~Bo Ning\footnote{Center for Applied Mathematics,
Tianjin University, Tianjin, 300072, P.R. China. Email: bo.ning@tju.edu.cn. Partially supported by NSFC grant 11601379.}}

\maketitle

\begin{abstract}
In this paper we determine the chromatic number
of graphs with two odd cycle lengths.
Let $G$ be a graph and $L(G)$ be the set of
all odd cycle lengths of $G$. We prove that: (1)
If $L(G)=\{3,3+2l\}$, where $l\geq 2$,
then $\chi(G)=\max\{3,\omega(G)\}$; (2) If $L(G)=\{k,k+2l\}$,
where $k\geq 5$ and $l\geq 1$, then $\chi(G)=3$.
These, together with the case $L(G)=\{3,5\}$
solved in \cite{W},
give a complete solution to the general problem addressed
in \cite{W,CS,KRS}. Our results also improve a
classical theorem of Gy\'{a}rf\'{a}s which asserts that
$\chi(G)\le 2|L(G)|+2$ for any graph $G$.

\end{abstract}



\section{Introduction}
Only simple graphs are considered. For a graph $G$,
let $\chi(G)$, $\omega(G)$, and $L(G)$ denote
the chromatic number of $G$, the size of maximum cliques in $G$, and the set of
all odd cycle lengths of $G$, respectively.
For notations not defined, we refer the reader to \cite{BM}.

The study of the relation between $\chi(G)$, $\omega(G)$ and $L(G)$ is
a fundamental area in graph theory and has been a subject of extensive research.
It is well-known that $\chi(G)\leq2$ if and only if $L(G)=\emptyset$.
A general upper bound for $\chi(G)$ in terms of the size of $L(G)$ was proposed by Bollob\'{a}s and Erd\H{o}s \cite{E90},
where they conjectured that $\chi(G)\leq 2|L(G)|+2$ for any $G$.
In \cite{G}, Gy\'{a}rf\'{a}s confirmed this by showing that if $|L(G)|=k\geq 1$,
then $\chi(G)\leq 2k+2$ with equality if and only if some block of $G$ is a $K_{2k+2}$.
If one considers the elements of $L(G)$, then often the value of $\chi(G)$ can be improved.
Indeed, in \cite{W} Wang proved that $\chi(G)=3$ if $L(G)=\{k\}$ for some $k\geq 5$.
Kaiser, Ruck\'{y} and \v{S}krekovski \cite{KRS} obtained a slight improvement that any proper 3-coloring of
an odd cycle of $G$ can be extended to a proper 3-coloring of $G$, assuming $G$
contains no $K_4$ and has $|L(G)|=1$.
The problem of determining $\chi(G)$ seems to be much harder for graphs with $|L(G)|=2$.
The case $L(G)=\{3,5\}$ was resolved by Wang \cite{W}, where he proved that if $G$ contains neither $K_4$ nor
$W_6$ (a wheel on six vertices) then $\chi(G)=3$, and otherwise $\chi(G)=\max\{4,\omega(G)\}$.
In \cite{CS}, Camacho and Schiermeyer showed that every graph $G$ with $L(G)=\{k,k+2\}$ for $k\geq 5$
satisfies $\chi(G)\leq 4$. The special case $L(G)=\{5,7\}$
was improved to $\chi(G)=3$ by Kaiser, Ruck\'{y} and \v{S}krekovski in \cite{KRS}.

In this paper, we determine $\chi(G)$ for every graph $G$ with $|L(G)|=2$.
Our main theorems are as follows.
\begin{theorem}\label{Thm:3l}
Let $l\geq 2$ be an integer. Any graph $G$ with $L(G)=\{3,3+2l\}$ has $\chi(G)=\max\{3,\omega(G)\}$.
\end{theorem}

\begin{theorem}\label{Thm:kl}
Let $k\geq 5$ and $l\geq 1$ be integers. Any graph $G$ with $L(G)=\{k,k+2l\}$ has $\chi(G)=3$.
\end{theorem}

We point out that these results improve the aforementioned theorem of Gy\'{a}rf\'{a}s in the family of graphs considered.
Recently, the theorem of Gy\'{a}rf\'{a}s was extended to cycles of consecutive odd lengths in a joint
paper \cite{LM} of the first author. Answering a conjecture of Erd\H{o}s \cite{E92},
Kostochka, Sudakov and Verstra\"ete in \cite{KSV} proved that every triangle-free graph $G$ with $|L(G)|=k$
satisfies $\chi(G)=O(\sqrt{k/\log k})$. For general $L(G)$, the precise value of $\chi(G)$ seems to be out of reach.
However, maybe it is possible to determine the maximum integer $t$ such that any triangle-free graph
$G$ with $|L(G)|=t$ has $\chi(G)=3$. The Gr\"otzsch graph and Chv\'atal graph both have $L(G)=\{5,7,9,11\}$ and $\chi(G)=4$,
which, together with Theorem \ref{Thm:kl}, show that $2\leq t\leq 3$.
It will be interesting to see if $t=3$.

Let $G=(V,E)$ be a graph, $x,y$ be vertices of $G$, and $H, H'$ be subgraphs of $G$.
For a subset $S$ of $V$, by $N_H(S)$ we denote the set of vertices in $V(H)\backslash S$,
each of which is adjacent to some vertex of $S$ in $G$.
We also denote by $H-S$ (and $H-H'$, respectively) the induced subgraph of $H$ on the vertex set $V(H)\backslash S$
(and $V(H)\backslash V(H')$, respectively).
For $x,y\in V(H)$, the \emph{distance} in $H$ between $x$ and $y$,
denoted by $d_{H}(x,y)$, is the length of a shortest path in $H$ with endpoints $x$ and $y$.
For a cycle or a path $Q$, the {\em length} of $Q$, denoted by $|Q|$, counts the number of edges in $Q$.
A cycle $C$ is called a $k$-\emph{cycle} if $|C|=k$.
If we draw a cycle $C$ as a circle in the plane, then $xCy$ denotes
the path on $C$ from $x$ to $y$ in the clockwise direction.
A path $P$ with endpoints $x$ and $y$ is called an {\em $(x,H,y)$-path}
if $V(P)\backslash\{x,y\}\subseteq V(H)$, and an {\em $(H,H')$-path} if $V(P\cap H)=\{x\}$ and $V(P\cap H')=\{y\}$.
For the convenience, we use $\widehat{P}$ to denote $P-\{x,y\}$.
An $H$-\emph{bridge} of $G$ is either an edge with two endpoints in $V(H)$
or a subgraph induced by a component $D$ of $G-H$ together with all edges between $D$ and $H$.
For subsets $A,B$ of $V$, the pair $(A,B)$ is called a \emph{$k$-separation} of $G$ if $A\cup B=V$,
$|A\cap B|=k$, and $G$ has no edges between $A\backslash B$ and $B\backslash A$.
A graph $G$ is $k$-{\em chromatic} if $\chi(G)=k$, and is $k$-\emph{critical} if $G$ is $k$-chromatic
but any proper subgraph of $G$ is not. If there is no danger of ambiguity, we often do not distinguish
the vertex set and the graph induced by it.
And if $H$ consists of a single vertex $v$, we also often write $v$ instead of $H$ or $\{v\}$ in the above notations.

The organization of this paper is as follows. In Section 2, we prove Theorem \ref{Thm:3l}.
In Section 3, we prove Theorem \ref{Thm:kl}, assuming Lemmas \ref{Lem:3connected} and \ref{Lem:2oddcycle}.
We then complete the proofs of Lemmas \ref{Lem:3connected} and \ref{Lem:2oddcycle} in Sections \ref{sec:3connected} and \ref{sec:2oddcycle}, respectively.

\section{Proof of Theorem \ref{Thm:3l}}
Throughout this section, let $G$ be a graph with
\begin{equation}
\label{equ:3l}
L(G)=\{3,k\}, \text{ where } k:=3+2l \text{ and } l\ge 2.
\end{equation}
We shall show that $\chi(G)=\max\{3,\omega(G)\}$.
It is fair to assume that $G$ is 2-connected.
Otherwise, there is
a cut vertex $u$ such that $G_1\cup G_2=G$ and $G_1\cap G_2=\{u\}$. Assume
that $k\in L(G_1)$. Then $L(G_1)$ is either $\{k\}$ or $\{3,k\}$,
and $L(G_2)$ can be $\emptyset$, or $\{3\}$, or $\{k\}$, or $\{3,k\}$.
Then we can use induction for $L(G_i)=\{3,k\}$, or Wang's result \cite{W} that $\chi(G_i)=3$ for $L(G_i)=\{k\}$,
or Gy\'arf\'as' result \cite{G} that $\chi(G_2)=\max\{3,\omega(G_2)\}$ for $L(G_2)=\{3\}$.
Putting the above together, it will be easy to see that $\chi(G)=\max\{3,\omega(G)\}$.

By \eqref{equ:3l}, observe that $\omega(G)\in \{3,4\}$. According to the value of $\omega(G)$,
we divide the proof of Theorem \ref{Thm:3l} into two subsections as follows.

\subsection{$\omega(G)=4$}
Let $X$ be a $K_4$ in $G$ with $V(X)=\{x_1,x_2,x_3,x_4\}$.
We will need to prove $\chi(G)=4$.
To achieve this, we propose to show that for any component $H$ in $G-X$,
any proper 4-coloring of $X$ can be extended to a proper 4-coloring of $G[V(X\cup H)]$.

First we claim that for distinct $x_i,x_j\in V(X)$ there is no $(x_i,x_j)$-path of even length in $G$ internally disjoint from $X$.
Suppose to the contrary that there is a such path $P$ in $G$, say from $x_1$ to $x_2$. Then $P\cup x_1x_2$
and $P\cup x_1x_3x_4x_2$ are two odd cycles in $G$ with lengths differ by two, a contradiction to \eqref{equ:3l}.
This proves the claim.

Suppose that $H$ contains an odd cycle, say $C$. Since $G$ is 2-connected, there are two disjoint
$(X,C)$-paths $P_1,P_2$, say from $x_1,x_2\in V(X)$ to $y_1,y_2\in V(C)$, respectively.
Since $|C|$ is odd, there exists a $(y_1,y_2)$-path $Q$ on $C$ such that
$P_1\cup P_2\cup Q$ is an even $(x_1,H,x_2)$-path in $G$, a contradiction.
So, $H$ is bipartite.

Let $(A,B)$ be the bipartition of $H$.
Next we show that no distinct $x_i,x_j\in V(X)$ can be adjacent to the same part in $(A,B)$.
Otherwise, by symmetry we may assume that there exist $a\in A\cap N_G(x_1)$ and $a'\in A\cap N_G(x_2)$.
Let $P$ be an $(a,a')$-path of $H$. As $|P|$ is even, we see $x_1a\cup P\cup a'x_2$ is
an even $(x_1,H,x_2)$-path in $G$, a contradiction to the claim.

We can then derive that there are at most two vertices in $X$ adjacent to $H$,
say $V(X)\cap N_G(H)\subseteq \{x_1,x_2\}$.
Now it is clear that any proper 4-coloring $\varphi$ of $X$ can be extended to a proper 4-coloring of $G[V(X\cup H)]$,
by coloring all vertices of $A$ by the color $\varphi(x_3)$ and all vertices of $B$ by the color $\varphi(x_4)$.
The proof of Theorem \ref{Thm:3l} when $\omega(G)=4$ is completed. {\hfill$\Box$}

\subsection{$\omega(G)=3$}

To finish the proof of Theorem \ref{Thm:3l}, it remains to consider a graph $G$ containing no $K_4$.
We are going to prove $\chi(G)=3$ by the means of contradiction.
Let $G$ be a minimal $K_4$-free graph satisfying \eqref{equ:3l} but $\chi(G)\ge 4$. We claim that $G$ is 4-critical.
Indeed, if not, then there exists $e\in E(G)$ such that $\chi(G-e)\geq 4$;
by the choice of the minimality of $G$, we have $L(G-e)=\{3\}$ or $\{k\}$,
which, by Gy\'arf\'as' result \cite{G} or Wang's result \cite{W}, implies that $\chi(G-e)=3$, a contradiction.
So
\begin{equation*}
G \text{ is 4-critical, which implies that } \delta(G)\geq 3 \text{ and } G \text{ is 2-connected.}
\end{equation*}
Recall that we write $k=3+2l$, where $k\ge 7$ (as $l\ge 2$).

Our starting point is a result of Voss \cite[Theorem 2]{Voss82} (also see \cite{Voss}) that
every $K_4$-free graph with chromatic number at least 4 contains an odd cycle with at least two diagonals.
By this theorem, $G$ contains a $k$-cycle $C$ with at least two diagonals, as clearly such cycle can not be a triangle.
Let $C:=v_0v_1\ldots v_{k-1}v_0$, and $G_0:=G[V(C)]$. (The subscripts will be taken modulo $k$ in the rest of this section.)

In what follows, we will prove a sequence of claims.
The first claim shows that the induced subgraph $G_0$ consists of the $k$-cycle $C$ and exactly two diagonals.
Without loss of generality, we may assume that

\setcounter{claim}{0}
\begin{claim}
$E(G_0)=E(C)\cup \{v_0v_2,v_1v_3\}$.
\end{claim}

\begin{pf}
For any diagonal $v_iv_j$ of $C$, there exists a $(v_i,v_j)$-path $P$ on $C$ such that
$P\cup v_iv_j$ forms an odd cycle. Since $L(G)=\{3,k\}$ and $j\notin \{i-1,i+1\}$,
we see that $P\cup v_iv_j$ is of length less than $k$ and thus a triangle, implying that $j\in \{i-2,i+2\}$.
Without loss of generality let $v_0v_2$ be a diagonal of $C$.
Consider any other diagonal $v_iv_{i+2}$ of $C$.
If $i\notin \{1,k-1\}$, then there is a $(k-2)$-cycle $(C-\{v_1,v_{i+1}\})\cup v_0v_2\cup v_iv_{i+2}$,
so $l=1$, a contradiction. Thus, except $v_0v_2$, only $v_1v_3$ or $v_1v_{k-1}$ can be a diagonal of $C$,
and one can easily see that both of them cannot be. This proves Claim 1.
\end{pf}

\medskip

We define a proper 3-coloring $\varphi:V(G_0)\to \{1,2,3\}$ of $G_0$ by the following rule:
\begin{itemize}
\item Let $S_1:=\{v_3,v_5,...,v_{k-2},v_0\}$ and $S_2:=\{v_2,v_4,...,v_{k-1}\}$.
\item Assign $\varphi(v_1):=3$, and for any $j\in \{1,2\}$ and $x\in S_j$, assign $\varphi(x):=j$.
\end{itemize}
The essential idea behind the coming claims is to show that for every component $H$
in $G-V(G_0)$,
\begin{equation}
\label{equ:3lextension}
\varphi \text{ can be extended to a proper 3-coloring of } G[V(G_0\cup H)].
\end{equation}
Note that, if true, this in turn will give rise to a proper 3-coloring of $G$ and complete the proof of Theorem \ref{Thm:3l}.
We prove by contradiction. Suppose that there exists a component $H$ in $G-V(G_0)$ such that (2) does not hold.
\begin{claim}
If there exist $v_i,v_j\in N(u)\cap V(C)$ for some $u\in V(H)$, then $d_C(v_i,v_j)=1$ or $2$.
Moreover, $\{v_i,v_j\}\neq \{v_p,v_{p+2}\}$ for any $p\in \{k-1,0,1,2\}$.
\end{claim}

\begin{pf}
There exists a $(v_i,v_j)$-path $P$ on $C$ such that $P\cup v_iu\cup uv_j$ forms an odd cycle.
As $L(G)=\{3,k\}$, it is easy to see that $d_C(v_i,v_j)=1$ or $2$.
Suppose that $\{v_i,v_j\}=\{v_p,v_{p+2}\}$ for some $p\in \{k-1,0,1,2\}$.
Then it is easy to check that this will force a 5-cycle in $G$, a contradiction.
\end{pf}

\begin{claim}
We may assume that $|V(H)|\ge 2$.
\end{claim}
\begin{pf}
Suppose to the contrary that $V(H)=\{u\}$ for some $u\in V(G)$.
Claim 2 shows that any two neighbors of $u$ is of distance one or two on $C$.
Since $\delta(G)\ge 3$ and $|C|=k\ge 7$, one can deduce that $N(u)=\{v_i,v_{i+1},v_{i+2}\}$ for some $i$.
If $v_1\notin N(u)$, then we can assign $\varphi(u):=3$ such that \eqref{equ:3lextension} holds.
So $v_1\in N(u)$, which means that $i\in \{k-1,0,1\}$, contradicting Claim 2.
\end{pf}

\medskip

\begin{claim}
Let $v_i,v_j\in V(C)$ with $i\neq j$. If there are two $(v_i,H,v_j)$-paths $P$ and $Q$ with lengths differ by one,
then $\{v_i,v_j\}\cap \{v_1,v_2\}=\emptyset$, $\{|P|,|Q|\}=\{l+1,l+2\}$, and
the $(v_i,v_j)$-path on $C$ containing $\{v_1,v_2\}$ is of length $l+2$.
\end{claim}
\begin{pf}
Recall that $k=2l+3$. Let $p:=|P|$ and $q:=|Q|$, and assume by symmetry that $p$ is odd.
Then $p\ge 3$, implying that $p+q\ge 5$.
Let $X$ be the even $(v_i,v_j)$-path on $C$ such that $C_1:=X\cup P$ forms an odd cycle.
Then $C_2:=(C-\widehat{X})\cup Q$ also is an odd cycle. As $L(G)=\{3,k\}$,
$|C_1|+|C_2|=|C|+p+q\in \{6,k+3,2k\}$. In view of $p+q\ge 5$, we see that $|C_1|+|C_2|=2k$ and thus $p+q=k$.

We first show that $\{v_i,v_j\}\cap \{v_1,v_2\}=\emptyset$. Suppose not, say $v_i=v_1$.
If $v_j\in \{v_0,v_2\}$, then $(C-v_iv_j)\cup P$ is an odd cycle of length more than $k$, a contradiction.
If $v_j=v_3$, then $(C-\{v_1,v_2\})\cup v_0v_2v_1\cup P$ is an odd cycle of length more than $k$, a contradiction.
So $v_j\in V(C)-\{v_0,v_1,v_2,v_3\}$.
Let $X$ be a $(v_1,v_j)-$path on $C$ through $v_0$, and choose $R\in \{P,Q\}$
with the parity different from $X$. Then $X\cup R$ and $(X-v_0v_1)\cup R\cup v_1v_3v_2v_0$ are two odd cycles
whose lengths differ by two, a contradiction to (1).
This proves $\{v_i,v_j\}\cap \{v_1,v_2\}=\emptyset$.

Let $Z$ be the $(v_i,v_j)$-path on $C$ containing $\{v_1,v_2\}$. So $|Z|\ge 3$.
Let $\{R_1,R_2\}=\{P,Q\}$ such that $C':=R_1\cup Z$ forms an odd cycle.
Then $C'':=R_2\cup (Z-\{v_0v_1,v_1v_2\})\cup v_0v_2$ is also an odd cycle.
So $|C'|+|C''|=p+q+2|Z|-1=k+2|Z|-1\in \{6,k+3,2k\}$.
As $|Z|\ge 3$, this shows that $|Z|=(k+1)/2=l+2$ and $|C'|=|C''|=k$,
further implying that $|R_1|=l+1$ and $|R_2|=l+2$. Claim 4 is proved.
\end{pf}

\medskip

A \emph{book} of $r$ pages, denoted by $B^{*}_r$, is a graph consisting
of $r$ triangles sharing with one common edge. It was proved in \cite[Theorem~8]{W}
that every 2-connected non-bipartite graph containing no odd cycles other than 3-cycles is either a $K_4$ or a book.
This leads us to the next claim.

\begin{claim}
Every non-bipartite block in $H$ is a book $B^*_r$ for some $r\ge 1$.
\end{claim}

\begin{pf}
Let $B$ be a non-bipartite block in $H$. Suppose that $B$ contains a $k$-cycle $C'$, which is disjoint from $C$.
As $G$ is 2-connected, there exist two disjoint paths $X,Y$ from $x,y\in V(C)$ to $x',y'\in V(C')$,
respectively and internally disjoint from $C\cup C'$.
Let $P$ be an $(x,y)$-path on $C$ and $P'$ be an $(x',y')$-path on $C'$ such that $C_1:=P\cup P'\cup X\cup Y$ is an odd cycle.
Then $C_2:=(C-\widehat{P})\cup (C'-\widehat{P'})\cup X\cup Y$ is also odd. But $|C_1|+|C_2|=|C|+|C'|+2|X|+2|Y|>2k$, a contradiction to $L(G)=\{3,k\}$.
This shows that $B$ contains no $k$-cycles and thus $L(B)=\{3\}$.
Claim 5 then follows from Theorem 8 in \cite{W} just mentioned and the fact that $G$ is $K_4$-free.
\end{pf}

\begin{claim}
$H$ has at most one non-bipartite block.
\end{claim}

\begin{pf}
Suppose to the contrary that $H$ has two such blocks, say $B_1$ and $B_2$.
Let $W$ be a path in $H$ from $w_1\in V(B_1)$ to $w_2\in V(B_2)$ internally disjoint from $B_1\cup B_2$,
where $w_i$ is a cut-vertex of $H$ contained in $B_i$ for $i=1,2$.
By Claim 5, $B_i$ is a book and thus contains a triangle, say $T_i:=G[\{w_i,x_i,y_i\}]$.
Each of $x_i,y_i$ is either a vertex of degree two in the book $B_i$ or adjacent to such a vertex in $V(B_i)-V(T_i)$;
while for each vertex $u$ of degree two in the book $B_i$, there is a path from $u$ to $V(C)$ internally disjoint from $B_i$
(since $\delta\geq 3$).
Hence, by symmetry, we may assume that there exist two internally disjoint paths $P_1, P_2$ from $x_1,x_2$ to $v_i,v_j\in V(C)$
for some $i,j$, respectively and internally disjoin from $T_1\cup T_2\cup W\cup C$.

If one can choose the above $P_1,P_2$ such that $v_i\neq v_j$, then we can find three $(v_i,H,v_j)$-paths, namely,
$P:=P_1\cup x_1w_1\cup W\cup w_2x_2\cup P_2$, $(P-x_1w_1)\cup x_1y_1w_1$ and $(P-\{x_1w_1,x_2w_2\})\cup x_1y_1w_1\cup x_2y_2w_2$,
with three consecutive lengths, a contradiction to Claim 4. Thus, for all choices of $\{P_1,P_2\}$,
\begin{equation}
\label{equ:intersect}
P_1,P_2 \text{ intersect } V(C) \text{ at the same vertex, say } v_i.
\end{equation}
Then we get three cycles of consecutive lengths, which implies that the middle cycle is a $k$-cycle and so
\begin{equation}
\label{equ:midlength}
|P_1|+|P_2|+|W|+3=2l+3.
\end{equation}

Since $G$ is 2-connected, there exists a path $Q$ from $v_q\in V(C)-\{v_i\}$ to $w\in V(T_1\cup T_2\cup W\cup P_1\cup P_2)$,
internally disjoint from $C\cup T_1\cup T_2\cup W\cup P_1\cup P_2$.
If $w\in V(P_1\cup P_2)\cup \{x_1,y_1,x_2,y_2\}$, then we get a contradiction to \eqref{equ:intersect}.
Thus, $w\in V(W)$.
For each $i\in \{1,2\}$, let $Q_i:=Q\cup wWw_i\cup w_ix_i\cup P_i$,
then $Q_i$ and $(Q_i-w_ix_i)\cup w_iy_ix_i$ are two $(v_q,H,v_i)$-paths whose lengths differ by one.
Claim 4 then implies that for $i\in \{1,2\}$, the length of $Q_i$ is
$$|Q|+|wWw_i|+1+|P_i|=l+1.$$
Adding $|Q_1|$ and $|Q_2|$ up, we have $$|P_1|+|P_2|+|W|+2|Q|+2=2l+2,$$
which, compared with \eqref{equ:midlength}, shows that $|Q|=0$, a contradiction. This proves Claim 6.
\end{pf}

\medskip

By Claims 5 and 6, let $D$ be the unique non-bipartite block of $H$ (if existing), such that $V(D)=\{x_1,x_2,y_1,...,y_r\}$ and
$E(D)=\{x_1x_2\}\cup \{x_iy_j:1\le i\le 2, 1\le j\le r\}$, where $r\geq 1$.
Denote $H':=H-\{x_1x_2\}$ if $D$ exists; otherwise, denote $H':=H$. Therefore,
\begin{equation}
\label{equ:H'}
H' \text{ is connected and bipartite.}
\end{equation}
Let $(A,B)$ be the bipartition of $H'$. So $\{x_1,x_2\} \subseteq A$ or $\{x_1,x_2\} \subseteq B$ if $D$ exists.

\begin{claim}
If $N_H(v_1)\neq \emptyset$, then $D$ does not exist and thus $H=H'$ is bipartite.
\end{claim}

\begin{pf}
Suppose that $N_H(v_1)\neq \emptyset$ and there is the non-bipartite block $D$ of $H$.
Let $T$ be a triangle in $D$ and denote $V(T)$ by $\{x_1,x_2,x_3\}$.
Since $G$ is 2-connected, there exist two disjoint paths $P,Q$ from $V(C)$ to $V(T)$ internally disjoint from $C\cup T$.
Since $N_H(v_1)\neq \emptyset$, by rerouting paths if necessarily, we may assume that $P,Q$ are from $v_1,v_i\in V(C)$
to two vertices $u,v$ in $T$. Let $w \in V(T)-\{u,v\}$. Then $P \cup uv \cup Q$ and $P \cup uwv \cup Q$ are two $(v_1,H,v_i)$-paths with lengths differ by one,
however it is a contradiction to Claim 4.
\end{pf}

\medskip

Recall the sets $S_1=\{v_3,v_5,...,v_{k-2},v_0\}$ and $S_2=\{v_2,v_4,...,v_{k-1}\}$, and the proper 3-coloring $\varphi$ on $G_0$.

\begin{claim}
$N_H(v_1)=\emptyset$.
\end{claim}

\begin{pf}
We prove this claim by showing that if $N_H(v_1)\neq\emptyset$, then $\eqref{equ:3lextension}$ holds.
Without loss of generality, assume that there exists $u_1\in N_G(v_1)\cap A$. By Claim 7, $H$ is bipartite.

We first show that $N_H(S_1)\subseteq B$. Otherwise, there exists $v_iu_i\in E(G)$ for some $v_i\in S_1$ and $u_i\in A$.
Let $P$ be a $(v_1,H,v_i)$-path with even length. If $v_i=v_0$, then $P\cup v_0v_1$ and $P\cup v_1v_3v_2v_0$
are two odd cycles with lengths differ by two, a contradiction. So $v_i\in S_1-\{v_0\}$.
Let $X$ be the $(v_3,v_i)$-path on $C$ not containing $v_1$. Note that $X$ is even.
Then $v_1v_3\cup X\cup P$ and $v_1v_0v_2v_3\cup X\cup P$ are two odd cycles with lengths differ by two, which cannot be.

Next we show that $N_H(S_2)\subseteq A$. Suppose to the contrary that there exists $v_ju_j\in E(G)$ for some $v_j\in S_2$ and $u_j\in B$.
Let $Q$ be a $(v_1,H,v_j)$-path with odd length at least three.
If $v_j=v_2$, then $(C-v_1v_2)\cup Q$ is an odd cycle of length at least $k+2$, a contradiction.
So $v_j\in S_2-\{v_2\}$. Let $Y$ be the $(v_3,v_i)$-path on $C$ not containing $v_1$. So $Y$ is odd.
Then $v_1v_3\cup Y\cup Q$ and $v_1v_0v_2v_3\cup Y\cup Q$ are two odd cycles with lengths differ by two, again a contradiction.

Note that $V(C)=\{v_1\}\cup S_1\cup S_2$ and $G$ is 2-connected. So $N_H(S_1)\cup N_H(S_2)$ is not empty.
Recall that we have proved $N_H(S_1)\subseteq B$, $N_H(S_2)\subseteq A$ and $H$ is bipartite. So $\varphi$ can be extended to a proper 3-coloring of $G[V(G_0\cup H)]$,
by simply coloring all vertices in $A$ using color 1 and all vertices in $B$ using color 2.
\end{pf}

\begin{claim}
If there exist distinct $v_p,v_q\in S_i$ for some $i$ adjacent to $u_p\in A, u_q\in B$, respectively,
then the $(v_p,v_q)$-path on $C$ not containing $v_1$ is of length $l+1$, any $(u_p,u_q)$-path in $H'$ is of length $l$, $l$ is odd, and $|N_G(A)\cap S_i|=|N_G(B)\cap S_i|=1$.

\end{claim}

\begin{pf}
By \eqref{equ:H'}, any $(u_p,u_q)$-path $P$ in $H'$ is of odd length .
Let $X$ be the $(v_p,v_q)$-path on $C$ not containing $v_1$, and $Y$ be the $(v_p,v_q)$-path $(C-\widehat{X}-\{v_1\})\cup v_0v_2$.
By the definitions of $S_1$ and $S_2$, both $X$ and $Y$ are even with $|X|+|Y|=k-1$.
Then $X\cup v_pu_p\cup P\cup u_qv_q$ and $Y\cup v_pu_p\cup P\cup u_qv_q$ are two odd cycles,
implying that $|X|+|Y|+2|P|+4\in \{6,k+3,2k\}$.
As $|X|+|Y|=k-1$ and $|P|\ge 1$, we deduce that $|X|+|Y|+2|P|+4=2k$.
This implies that $|P|=(k-3)/2=l$ and $|X|=(k-1)/2=l+1$.

Suppose that there is some $v_j\in N_G(A)\cap S_i-\{v_p\}$.
Note that $|C|=2l+3$. By a similar argument above, we have $d_C(v_j,v_q)=d_C(v_p,v_q)=l+1$.
Since $v_j\neq v_p$, vertices $v_p,v_q,v_j$ must lie on $C$ in cyclic order
and thus the $(v_p,v_j)$-path on $C$ containing $v_q$ is of length $2l+2$, a contradiction to the definition of $S_i$.
This shows that $|N_G(A)\cap S_i|=1$ and similarly $|N_G(B)\cap S_i|=1$, completing the proof.
\end{pf}

\begin{claim}
$l$ is odd.
\end{claim}

\begin{pf}
Suppose for a contradiction that $l$ is even. By Claim 9, we see $N_H(S_i)\subseteq A$ or $B$ for each $i$.
By the symmetry between $A$ and $B$, we have two cases (see below) to consider;
and we will show that in each case, $\varphi$ can be extended to a proper 3-coloring of $G[V(G_0 \cup H)]$.
Note that $N_H(v_1)=\emptyset$ by Claim 8.

Suppose $N_H(S_1)\subseteq A$ and $N_H(S_2)\subseteq B$. We may further assume that $x_1,x_2\in A$ (if $D$ exists).
Then we can extend $\varphi$ onto $G[V(G_0\cup H)]$, by coloring $x_1$ using color 3, all vertices of $A-\{x_1\}$ using color 2 and all vertices of $B$ using color 1.
If $D$ does not exist, color each vertex in $A$ by 2 and each vertex in $B$ by 1.

Now we may assume $N_H(S_1)\cup N_H(S_2)\subseteq A$. Suppose that $D$ exists. If $x_1,x_2\in B$, we can color $x_1$ by color 1, color $B-\{x_1\}$ by color 2, and
color all vertices of $A$ by color 3. Thus, $x_1,x_2\in A$. If there exist some $i,j\in \{1,2\}$ such that $x_i\notin N_H(S_j)$,
then we can color $x_i$ by color $j$, color all vertices of $A-\{x_i\}$ by color $3$, and color all vertices of $B$ by color $3-j$.
It remains to consider the situation that for each $i\in \{1,2\}$, there exist $v_p\in S_1$ and $v_q\in S_2$ such that $v_p,v_q\in N(x_i)$.
By Claim 2, we see that $d_C(v_p,v_q)=1$ or $2$. As $\{v_p,v_q\}\neq \{v_0,v_2\}$ (by Claim 2) and $v_p\in S_1$, $v_q\in S_2$, it holds that in fact $d_C(v_p,v_q)=1$.
Hence, we may assume that there exist vertices $v_s,v_{s+1}\in N(x_1)$ and $v_t,v_{t+1}\in N(x_2)$.
Clearly $s \neq t$, for otherwise $G$ contains a $K_4$. Then $(C-\{v_sv_{s+1},v_tv_{t+1}\})\cup (v_sx_1v_{s+1})\cup (v_tx_2v_{t+1})$ forms a $(k+2)$-cycle in $G$, a contradiction.
Suppose that $D$ does not exist. Then we color all vertices in $A$ by 3 and color all vertices in $B$ by 1 or 2.
This proves Claim 10.
\end{pf}

\begin{claim}
If there are $v_i\in S_1, v_j\in S_2$ both adjacent to $F$ for some $F\in \{A,B\}$,
then $d_C(v_i,v_j)=1$ and $N_G(v_i)\cap F=N_G(v_j)\cap F=\{u\}$ for some vertex $u$.
Moreover, if $N_G(S_1) \cap F \neq \emptyset$ and  $N_G(S_2) \cap F\neq \emptyset$ for some
$F \in \{A,B\}$, then $N_G(S_1) \cap V(H) = N_G(S_2) \cap V(H)$ and $\lvert N_G(S_1) \cap V(H) \rvert=1$.
\end{claim}

\begin{pf}
Let $P$ be any path in $H'$ from $u_i\in N_G(v_i)\cap F$ to $u_j\in N_G(v_j)\cap F$.
Clearly $|P|$ is even. Let $X$ be the $(v_i,v_j)$-path on $C$ not containing $v_1$, and $Y$ be the $(v_i,v_j)$-path $(C-\widehat{X}-\{v_1\})\cup v_0v_2$.
Then both $X$ and $Y$ are odd with $|X|+|Y|=k-1$, thus $C_1:=X\cup v_iu_i\cup P\cup u_jv_j$ and $C_2:=Y\cup v_iu_i\cup P\cup u_jv_j$ are two odd cycles.
This shows that $|X|+|Y|+2|P|+4\in \{6,k+3,2k\}$, so $|P|=0$ or $|P|=(k-3)/2=l$.
The latter case contradicts Claim 10, as $|P|$ is even. Hence $|P|=0$, implying that
both $N_G(v_i)\cap F$ and $N_G(v_j)\cap F$ consist of a single vertex, say $u$.
If $|C_1|=|C_2|=3$, then $|X|=|Y|=1$ and $k=3$, a contradiction.
If $|C_1|=|C_2|=k$, then $|X|=|Y|=k-2$ and thus $|X|+|Y|=2k-4=k-1$, implying that $k=3$, again a contradiction.
So we have $\{|C_1|,|C_2|\}=\{3,k\}$, which implies $\{|X|,|Y|\}=\{1,k-2\}$.
If $|Y|=1$, then $\{v_i,v_j\}=\{v_0,v_2\}$, a contradiction to Claim 2.
So $|X|=1$. That is, $d_C(v_i,v_j)=1$.

Suppose there is a vertex $w\in V(H)-\{u\}$ such that $w\in N_H(v_p)\cap N_H(v_{p+1})$ for some $p$.
If $\{v_i,v_{i+1}\}\neq \{v_p,v_{p+1}\}$, then $(C-\{v_iv_{i+1},v_pv_{p+1}\})\cup (v_iuv_{i+1})\cup (v_pwv_{p+1})$
is a $(k+2)$-cycle in $G$, a contradiction. So $\{v_i,v_{i+1}\}=\{v_p,v_{p+1}\}$.
Let $P'$ be a path in $H'$ from $u$ to $w$, and $Q:=v_iu\cup P'\cup wv_{i+1}$ be from $v_i$ to $v_{i+1}$ with $|Q|\ge 3$.
Denote $C'$ to be the cycle $(C-\{v_iv_{i+1}\})\cup Q$ if $Q$ is odd and $(C-\{v_iv_{i+1},v_0v_1,v_1v_2\})\cup v_0v_2\cup Q$ if $Q$ is even. As $|Q|\ge 3$, in either case $C'$ is an odd cycle of length more than $k$, a contradiction.
This proves Claim 11.
\end{pf}

\medskip

By Claim 11, if $N_G(S_1) \cap F \neq \emptyset$ and  $N_G(S_2) \cap F\neq\emptyset$ for some $F \in \{A,B\}$, then $N_G(S_1) \cap V(H)=N_G(S_2) \cap V(H)=\{u\}$ for some vertex $u$, and there exists a unique number $p \in \{2,3,...,k\}$ such that $u \in N_H(v_p) \cap N_H(v_{p+1})$ and we let $U:=\{u\}$; otherwise, let $U:=\emptyset$.

\begin{claim}
$N_H(S_i)\subseteq A\cup U$ and $N_H(S_j)\subseteq B\cup U$ for some $\{i,j\}=\{1,2\}$.
\end{claim}
\begin{pf}
By Claim 11, we see that if $N_H(S_i)\backslash U\neq \emptyset$ for each $i\in \{1,2\}$, then this assertion follows.
Thus, without loss of generality, we assume that $N_H(S_2)\subseteq U$.
If $N_H(S_1)\subseteq A\cup U$ or $B\cup U$, then again this assertion follows.
So we may assume that there exist $v_i\in N_G(A-U)\cap S_1$ and $v_j\in N_G(B-U)\cap S_1$
(and we will choose distinct $v_i,v_j$ if existing).
Suppose $U=\{u\}$. Recall $u$ is adjacent to both $v_p$ and $v_{p+1}$.
By the symmetry, assume that $u\in A$ and $v_p\in S_2$.
We then apply Claim 11 to the pair of vertices $v_i,v_p$, and it follows that $N_G(v_i)\cap A=\{u\}$, a contradiction to the choice of $v_i$. So $U=\emptyset$. Then we have $(S_2\cup \{v_1\})\cap N_G(H)=\emptyset$ and thus
$|S_1\cap N_G(H)|\ge 2$ (as $G$ is 2-connected), implying that $v_i\neq v_j$ (by the choice).
This in turns enables us to apply Claim 9 and conclude that $N_G(A)\cap S_1=\{v_i\}$ and $N_G(B)\cap S_1=\{v_j\}$.

By symmetry, if $D$ exists, then we assume $\{x_1,x_2\} \subseteq A$.
If $v_i$ is not adjacent to some vertex in $\{x_1,x_2\}$, say $x_1$,
then $\varphi$ can be extended onto $V(H)$ by coloring $x_1$ using color 1, all vertices in $A-\{x_1\}$ using color 2 and all vertices in $B$ using color 3. It is clear that $\eqref{equ:3lextension}$ holds.
So $v_i$ is adjacent to both $x_1$ and $x_2$. Since $H'$ is connected,
there exists a path $P$ in $H'$ from $w\in N_G(v_j)\cap B$ to some vertex in $\{x_1,x_2\}$, say $x_1$.
By Claim 9, $|P|=l$.
Then $v_jw\cup P\cup x_1v_i$ and $v_jw\cup P\cup (x_1x_2v_i)$ are two $(v_i,H,v_j)$-paths of lengths $l+2$ and $l+3$, respectively, a contradiction to Claim 4.
If $D$ does not exist, then color all vertices in $A$ by 2 and all vertices in $B$ by 3. This proves Claim 12.
\end{pf}

\medskip

Let $(i,j)=(1,2)$ in Claim 12. Note that $N_H(v_1)=\emptyset$ (by Claim 8).
We show how to extend $\varphi$ onto $V(H)$ and make $\eqref{equ:3lextension}$ hold.
By symmetry, if $x_1,x_2$ exist, then we assume $\{x_1,x_2\} \subseteq A$.
Suppose that either $U=\emptyset$, or $x_1,x_2$ do not exist, or $U=\{u\}$, $x_1,x_2$ exist and $x_ru \not \in E(G)$ for some $r \in \{1,2\}$.
Then we can color vertices in $\{x_r,u\}$ using color $3$, all vertices in $A-\{x_r,u\}$ using color 2 and all vertices in $B-\{x_r,u\}$ using color 1.
Hence, we may assume that vertices $x_1,x_2, u$ exist and induce a triangle in $H$.
As $G$ is 2-connected, there is a path $P$ in $G-\{u\}$ from some vertex $v_s$ in $V(C)$ to some vertex $x_t$ in $\{x_1,x_2\}$ internally disjoint from $V(C)$. By symmetry, we assume $x_r=x_1$.
Recall that $u$ is adjacent to both $v_p$ and $v_{p+1}$.
By the symmetry between $v_p$ and $v_{p+1}$, let $v_s\neq v_p$.
Then $v_pux_1\cup P$ and $v_pux_2x_1\cup P$ are two $(v_p,H,v_s)$-paths with lengths differ by one.
By Claim 4, $\{v_p,v_s\}\cap \{v_1,v_2\}=\emptyset$ and the $(v_p,v_s)$-path $X$ on $C$ containing $\{v_1,v_2\}$ is of length $l+2$.
This also shows $v_t\notin \{v_p,v_{p+1}\}$.
Then $v_{p+1}ux_1\cup P$ and $v_{p+1}ux_2x_1\cup P$ are two $(v_{p+1},H,v_s)$-paths with lengths differ by one as well.
By Claim 4 again, the $(v_{p+1},v_s)$-path on $C$ containing $\{v_1,v_2\}$ is of length $l+2$, which is a contradiction to $|X|=l+2$.
The proof of Theorem \ref{Thm:3l} is finished.
{\hfill$\Box$}

\section{Proof of Theorem \ref{Thm:kl}}
In this section, we shall prove Theorem \ref{Thm:kl}, assuming the following two lemmas whose proofs are postponed to the later sections.

\begin{lemma}\label{Lem:3connected}
Let $G$ be a 4-critical graph with $L(G)=\{k,k+2l\}$, where $k\geq 5$ and $l\geq 1$.
Then $G$ is 3-connected.
\end{lemma}

\begin{lemma}\label{Lem:2oddcycle}
Let $G$ be a 4-critical graph with $L(G)=\{k,k+2l\}$, where $k\geq 5$ and $l\geq 1$.
Then every two odd cycles in $G$ intersect in at least two vertices.
\end{lemma}

Like in the proof of Theorem \ref{Thm:3l}, we start the arguments by finding a cycle with certain property.
We say a cycle $C$ in $G$ is \emph{non-separating} if $G-V(C)$ is connected.
The coming result will be needed in the proof.
\begin{theorem}[\cite{TT,BV}]
\label{Thm:non-separating}
Every 3-connected non-bipartite graph contains a non-separating induced odd cycle.
\end{theorem}

Now we are ready to prove Theorem \ref{Thm:kl}.

\medskip

\noindent{\bf Proof of Theorem \ref{Thm:kl}.(Assuming Lemmas \ref{Lem:3connected} and \ref{Lem:2oddcycle})}
We prove by contradiction. Suppose it is not true. Then there exists a
counterexample graph $G$ such that the number of vertices is minimal,
and subject to this, the number of edges is minimal. So, similar as the proof of Theorem \ref{Thm:3l}, it is
\begin{equation*}
\text{ 4-critical and clearly non-bipartite with } L(G)=\{k,k+2l\}, \text{ where } k\ge 5 \text{ and } l\ge 1.
\end{equation*}
By Lemma \ref{Lem:3connected}, $G$ is 3-connected. Then by Theorem \ref{Thm:non-separating},
$G$ has a non-separating induced odd cycle $C$ such that $H:=G-V(C)$ is connected.
Moreover, Lemma \ref{Lem:2oddcycle} implies that $H$ is bipartite.
Let $(A,B)$ be the bipartition of $H$.
Since $\delta(G)\ge 3$,
\begin{equation}\label{equ:Cneighbor}
\text{every vertex on $C$ has at least one neighbor in } H.
\end{equation}
We will need to prove a sequence of claims and then arrive at the final contradiction to conclude this proof.

\setcounter{claim}{0}
\begin{claim}\label{CLNeig1}
For any $u\in V(C)$, $N_H(u)\subseteq A$ or $B$.
\end{claim}

\begin{pf}
Suppose that some $u\in V(C)$ has two neighbors $a\in A$ and $b\in B$.
Since $H$ is connected and bipartite,
there is an $(a,H,b)$-path of odd length. So $D:=ua\cup P\cup bu$ is
an odd cycle such that $V(C\cap D)=\{u\}$, contradicting Lemma \ref{Lem:2oddcycle}.
\end{pf}

\medskip

We can further deduce that

\begin{claim}\label{CLNeig2}
$N_{H}(C)\subseteq A$ or $B$.
In the rest of this proof, assume that $N_{H}(C)\subseteq A$.
\end{claim}

\begin{pf}
We say a vertex $u\in V(C)$ is of {\em type 0} if $N_{H}(u)\subseteq A$ and of {\em type 1} if $N_{H}(u)\subseteq B$.
In view of Claim \ref{CLNeig1}, every vertex on $C$ has type 0 or 1.

Suppose there exist vertices on $C$ of different types.
Then we can divide $C$ into paths $P_1,P_2,...,P_{2s}$ (appearing along a given cyclic order of $C$)
such that $V(C)=\bigcup_{i=1}^{2s}V(P_i)$
and for each $j\in \{0,1\}$, $V(P_{2i-j})$ consists of vertices of type $j$, where $1\leq i\leq s$.
We now define a 3-coloring $\varphi: V(G)\to \{0,1,2\}$ as follows:
every vertex in $A$ is colored by 1; every vertex in $B$ is colored by 0; and for every $j\in \{0,1\}$, we alternatively color $V(P_{2i-j})$ using colors $j,2$
such that the first vertex of the path (along the given cyclic order of $C$) is colored by $j$.
It is easy to see that $\varphi$ is a proper 3-coloring of $G$. This proves Claim \ref{CLNeig2}.
\end{pf}

\begin{claim}
$|C|=k$. Denote by $C:=x_0x_1x_2\cdots x_{k-1}x_0$.
\end{claim}

\begin{pf}
Suppose that $|C|=k+2l$. Write $C=x_0x_1x_2\cdots x_{k+2l-1}x_0$.
We show that for any $i$, $N_H(x_i)=N_H(x_{i+2})$. Otherwise, there
exist $a_1,a_2\in A$ with $a_1x_i,a_2x_{i+2}\in E(G)$. There is an $(a_1,a_2)$-path
$P$ with an even length in $H$. Thus $(C-\{x_{i+1}\})\cup \{a_1x_i,a_2x_{i+2}\}\cup P$ is
an odd cycle of length at least $k+2l+2$, a contradiction.
Now we can infer that in fact all $N_H(x_i)$ are the same set, implying that
there are triangles in $G$, a contradiction.
\end{pf}

\begin{claim}
$|V(H)|\geq 3$.
\end{claim}

\begin{pf}
Otherwise, $|V(H)|=1$ or $2$. Then by Claim \ref{CLNeig2}, in either case there exists a vertex $u\in V(H)$
which is adjacent to every vertex on $C$. This implies that there exist triangles in $G$, a contradiction.
\end{pf}

\begin{claim}
(1) If there is a vertex $y\in V(H)$ such that $x_iy,x_{i+2}y\in E(G)$, then $l=1$.
(2) If there is a trivial end-block (i.e., an edge) in $H$, then $l=1$.
\end{claim}
\begin{pf}
Set $x_j:=x_{i+2}$. Clearly $x_{i+1}y$, $x_{j+1}y\notin E(G)$,
since otherwise there is a triangle. So $x_{j+1}$ has a neighbor
$y'\in A-\{y\}$. There is an even $(y,y')$-path
$P$ in $H$, so $P\cup yx_{j}x_{j+1}y'$ and $P\cup yx_ix_{i+1}x_{j}x_{j+1}y'$
are two odd cycles with lengths differ by two. This proves (1).

Suppose $B:=yb$ is a trivial end-block in $H$, where $b$ is the cut-vertex. Since
$G$ is 3-connected, $y$ has two neighbors $x_i,x_j\in V(C)$. Since $|C|=k$
is the least odd cycle length, we have $d_C(x_i,x_j)=2$.
By Claim 5(1), we obtain $l=1$. This proves (2).
\end{pf}

\begin{claim}
For any 2-connected end-block $D$ in $H$, if there are two vertices
$x_i,x_{i+2}\in V(C)$ adjacent to $D$, then $l=1$.
\end{claim}

\begin{pf}
Suppose not. By Claim 5(1), assume that $l\geq 2$ and there are distinct $y_i,y_{i+2}\in V(D)$ such that
$x_iy_i,x_{i+2}y_{i+2}\in E(G)$.
Let $R$ be any $(y_i,y_{i+2})$-path in $H$, which must be of length $2l$.
This is because $R\cup \{x_iy_i,x_{i+2}y_{i+2}\}\cup (C-\{x_{i+1}\})$ is an odd
cycle of length $|R|+k=k+2l$.

Since $D$ is 2-connected, there are two disjoint $(y_i,y_{i+2})$-paths $P,Q$ in $D$
such that $|P|=|Q|=2l$, Then $C':=P\cup Q$ is an even cycle of
length $4l$. Write $C':=u_0u_1u_2\ldots u_{4l-1}u_0$
with $u_0:=y_i$ and $u_{2l}:=y_{i+2}$. Let $P_1:=x_iu_0$ and $P_2:=x_{i+2}u_{2l}$.
As $G$ is 3-connected, there exists a path $P_3$ from $v\in V(C)$ to $u_j\in V(C')-\{u_0,u_{2l}\}$,
internally disjoint from $P_1\cup P_2\cup C\cup C'$.
Next we aim to show
\begin{equation}
\label{equ:P3}
\text{ for every path } P_3 \text{ defined as above, } v=x_{i+1} \text{ and } u_j\in \{u_l,u_{3l}\}.
\end{equation}
By symmetry, assume that $0<j<2l$.
We draw $C'$ in the plane such that $u_0,u_1,...,u_{2l-1}$ appear on $C'$ clockwise,
and let $Q_1:=u_0C'u_j$, $Q_2:=u_jC'u_{2l}$ and $Q_3:=u_{2l}C'u_0$ so that $C'=Q_1\cup Q_2\cup Q_3$.

To prove \eqref{equ:P3}, we first show $v\in V(C)-\{x_i,x_{i+2}\}$. Otherwise, say $v=x_i$, then either
$C_1:=x_ix_{i+1}x_{i+2}\cup P_2\cup Q_2\cup P_3$ or $C_2:=(C-\{x_{i+1}\})\cup P_2\cup Q_2\cup P_3$ is odd. If $C_1$ is odd,
then $C_3:=P_1\cup Q_3\cup Q_2\cup P_3$ is also odd with $|C_3|-|C_1|=2l-2\in \{0,2l\}$, implying that $l=1$;
otherwise $C_2$ is odd, then $C_4:=(C-\{x_{i+1}\})\cup P_2\cup Q_3\cup Q_1\cup P_3$ is an odd cycle of length
at least $k+2l+1$, a contradiction. Now we see $P_1,P_2,P_3$ are disjoint paths.
Since $C$ is odd and $|Q_1|+|Q_2|=2l=|Q_3|$, there is a $(v,x_i)$-path $L$ on $C$
such that $C_5:=L\cup P_1\cup Q_1\cup P_3$, and $C_6:=L\cup P_1\cup (Q_2\cup Q_3)\cup P_3$
are odd. So $|C_6|-|C_5|=|Q_2|+|Q_3|-|Q_1|=4l-2|Q_1|\in \{0,2l\}$. Since $|Q_1|<2l$, this implies that
$|Q_1|=l$ and thus $u_j=u_l$.
Lastly, suppose that $v\neq x_{i+1}$, i.e., $v\in V(C)-\{x_i,x_{i+1},x_{i+2}\}$.
By the symmetry between $x_i$ and $x_{i+2}$, let $X$ be the $(v,x_{i+2})$-path on $C$ not containing $x_i$
such that $C_7:=X\cup P_2\cup Q_2\cup P_3$ is odd. Then
$C_8:=X\cup (x_{i+2}x_{i+1}x_i)\cup P_1\cup Q_1\cup P_3$ is also odd with $|C_8|-|C_7|=2$, implying $l=1$. This proves \eqref{equ:P3}.

Let $u_j=u_{l}$ and $Q_i$'s be as above.
Since $l\geq 2$, $Q_1-\{u_0,u_l\}$ is not empty. Since $G$ is 3-connected, there is a path
$R$ from $r\in V(Q_1)-\{u_0,u_l\}$ to $s\in (C'-Q_1)\cup C\cup P_1\cup P_2\cup (P_3-\{u_l\})$,
internally disjoint from $C'\cup C\cup P_1\cup P_2\cup P_3$.

If $s\in Q_2-\{u_l\}$, then $C'':=u_0Q_1r\cup R\cup sQ_2u_{2l}\cup Q_3$ is a cycle of length $4l$,
however the path $rQ_1u_l\cup P_3$ from $C''$ to $C$ contradicts \eqref{equ:P3}.
If $s\in Q_3-\{u_0,u_{2l}\}$, then $R_1:=u_0Q_1r\cup R\cup sQ_3u_{2l}$,
$R_2:=(C'-R_1)\cup R$ are two $(y_i,y_{i+2})$-paths in $H$, implying that $4l=|R_1|+|R_2|=|C'|+2|R|>4l$, a contradiction.
Hence, $s\notin C'$.
By \eqref{equ:P3}, we also have $s\notin (C-\{x_i,x_{i+2}\})\cup (P_3-\{u_l\})$. Therefore,
$s\in \{x_i,x_{i+2}\}$. In either case, let $C_1:=x_{i+1}s\cup R\cup rQ_1u_l\cup P_3$
and $C_2:=(C-\{x_{i+1}s\})\cup R\cup rQ_1u_l\cup P_3$. There is some
$C_i$, which is odd. As $C'$ is even,
the cycle $C'_i:=C_i\Delta C'$ is also odd. But $|C'_i|-|C_i|=
|C'|-2|rQ_1u_l|>4l-2l=2l$, a contradiction.
The proof of Claim 6 is completed.
\end{pf}

\begin{claim}
$l=1$ and thus $L(G)=\{k,k+2\}$, where $k\ge 5$.
\end{claim}

\begin{pf}
Suppose to the contrary that $l\ge 2$. By Claims 5 and 6,
we see that $H$ is not 2-connected, and all its end-blocks are 2-connected.
Let $D_1$ be an end-block of $H$,
$b\in V(D_1)$ be the cut-vertex of $H$ contained in $D_1$, and $D_2:=H-V(D_1-b)$.
Since $G$ is 3-connected, there exist $x_i\in V(C)$ and
$y_i\in V(D_1-b)$ such that $x_iy_i\in E(G)$.
Let $y_{i-1}, y_{i+1}\in V(H)$ such that $x_{i-1}y_{i-1}, x_{i+1}y_{i+1}\in E(G)$.
By Claim 5(1), $y_{i-1},y_{i+1}$ are distinct, and by Claim 6, $\{y_{i-1},y_{i+1}\}\not\subseteq V(D_1)$.
According to the locations of $y_{i-1}$ and $y_{i+1}$, we consider the following two cases.

Suppose that exactly one of $\{y_{i-1},y_{i+1}\}$ is in $D_2-b$,
say $y_{i-1}\in V(D_1)$ and $y_{i+1}\in V(D_2-b)$. Since $G$ contains no triangles, $y_{i-1}\neq y_i$.
Choose $y_{i-2},y_{i+2}\in V(H)$ such that $x_{i-2}y_{i-2},x_{i+2}y_{i+2}\in E(G)$.
We see that $y_{i-2},y_{i+2}\in V(D_2-b)$ (by Claim 6) and are distinct (as, otherwise, $G$ has an odd cycle of length $k-2$).
Let $P$ be a $(y_i,b)$-path in $D_1$,
$P_1$ a $(y_{i-2},b)$-path in $D_2$, and $P_2$ a $(y_{i+2},b)$-path in $D_2$.
Then by Claim \ref{CLNeig2}, $C_1:=(C-\{x_{i-1}\})\cup x_{i-2}y_{i-2} \cup P_1\cup P\cup y_ix_i$ and
$C_2:=(C-\{x_{i+1}\})\cup x_{i+2}y_{i+2}\cup P_2\cup P\cup y_ix_i$ are two odd cycles with
$|P_2|-|P_1|=|C_2|-|C_1|\in \{-2l,0,2l\}.$
Let $P'$ be a $(y_{i-1},b)$-path in $D_1$. Then
$C_3:=(y_{i-1}x_{i-1}x_{i-2}y_{i-2})\cup P_1\cup P'$ and
$C_4:=(y_{i-1}x_{i-1}x_ix_{i+1}x_{i+2}y_{i+2})\cup P_2\cup P'$
are two odd cycles satisfying that
$$|C_4|-|C_3|=2+|P_2|-|P_1|\in \{2-2l,2,2+2l\}\cap \{-2l,0,2l\}.$$
From the non-empty intersection, one can infer that $l=1$.

Now assume that $y_{i-1},y_{i+1}\in V(D_2-b)$. Let $Q$ be a $(y_i,b)$-path in $D_1$,
$Q_1$ a $(y_{i-1},b)$-path in $D_2$, and $Q_2$ a $(y_{i+1},b)$-path in $D_2$.
Consider the odd cycles $C_5:=(y_{i-1}x_{i-1}x_iy_i)\cup Q\cup Q_1$
and $C_6:=(y_{i+1}x_{i+1}x_iy_{i})\cup Q\cup Q_2$.
We can deduce that $|Q_2|-|Q_1|=|C_6|-|C_5|\in \{-2l,0,2l\}.$
Since $G$ is 3-connected, there exist $y_j\in V(D_1-b)$ and $x_j\in V(C)-\{x_i\}$ such that $x_jy_j\in E(G)$.
Let $Q_3$ be a $(y_j,b)$-path in $D_1$.
If $x_j$ is one of $\{x_{i-1},x_{i+1}\}$, then we are in the previous case.
So $x_j$ is distinct from $x_{i-1},x_{i+1}$.
Let $X$ be an $(x_j,x_{i-1})$-path on $C$ and $X'$ be an $(x_j,x_{i+1})$-path on $C$ such that both $|X|,|X'|$
are odd. By symmetry, let $|X'|-|X|=2$. Then $C_7:=X\cup x_jy_j\cup Q_3\cup Q_1\cup y_{i-1}x_{i-1}$,
and $C_8:=X'\cup x_jy_j\cup Q_3\cup Q_2\cup y_{i+1}x_{i+1}$ are two odd cycles with
$$|C_8|-|C_7|=|Q_2|-|Q_1|+(|X'|-|X|)\in \{2-2l,2,2+2l\}\cap \{-2l,0,2l\},$$
which again implies that $l=1$. This proves Claim 7.
\end{pf}

\medskip

In \cite{KRS} (see its Theorem~1.2), it was proved that every graph with $L=\{5,7\}$ has chromatic number 3.
By this result, we can assume that $k\geq 7$ in the rest of this section.

\begin{claim}
$H$ is not 2-connected.
\end{claim}

\begin{pf}
Suppose that $H$ is 2-connected. Note $C$ is the least odd cycle and $\delta(G)\geq 3$.
For any two consecutive vertices $x_i,x_{i+1}\in V(C)$, there are distinct $y_i,y_{i+1}\in A$ such
that $x_iy_i,x_{i+1}y_{i+1}\in E(G)$. There are 2 disjoint $(y_i,y_{i+1})$-paths $P_1,P_2$ in $H$, which are even.
Then for each $i=1,2$, $C_i:=P_i\cup (y_ix_ix_{i+1}y_{i+1})$ is an odd cycle, implying that $|P_i|\ge k-3$.
Then $C':=P_1\cup P_2$ forms an even cycle of length at least $2(k-3)\geq 8$, as $k\geq 7$.
Since $G$ is 3-connected, there are
3 disjoint paths $X_j$, $j=1,2,3$, from $u_j\in V(C')$ to $v_j\in V(C)$, internally disjoint
with $C\cup C'$. Let $C'_i$ be the $(u_{i-1},u_{i-2})$-path of $C'$, containing no $u_i$, where subscripts are taken mod 3.
Assume that $|C'_1|\geq |C'_2|\geq |C'_3|$. So $|C'_1|+|C'_2|-|C'_3|=|C'|-2|C'_3|\geq |C'|-2\lfloor\frac{|C'|}{3}\rfloor\geq \lceil\frac{|C'|}{3}\rceil\geq \lceil\frac{8}{3}\rceil=3$.

Since $C$ is odd and $C'$ is even, there exists a $(v_1,v_2)$-path $P$ on $C$ such that
$C_3:=P\cup X_1\cup X_2\cup C'_3$ and $C_4:=P\cup X_1\cup X_2\cup (C'_1\cup C'_2)$
are both odd. However, $|C_4|-|C_3|=|C'_1|+|C'_2|-|C'_3|\geq 3$,
contradicting $L(G)=\{k,k+2\}$. This proves Claim 8.
\end{pf}

\medskip

Let $x$ be a cut-vertex with $V(H_1\cap H_2)=\{x\}$ and $H_1\cup H_2=H$.
For a pair of vertices $\{x_i,x_{i+2}\}$ on $C$, we say that it is
\emph{feasible} (with respect to the cut-vertex $x$), if $N(x_i)\cap V(H_1-x)\neq \emptyset$,
and $N(x_{i+2})\cap V(H_2-x)\neq \emptyset$.

\begin{claim}
For any cut-vertex $x$ of $H$, $N(x)\cap V(C)=\emptyset$ and there exists a feasible pair $\{x_i,x_{i+2}\}$.
\end{claim}

\begin{pf}
If there exist $u,v\in N(x)\cap V(C)$,
then $u,v$ are of distance 2 on $C$, since otherwise there is an odd cycle of
length less than $k$. This shows that $|N(x)\cap V(C)|\leq 2$.
Assume, if existing, that $x_0,x_2\in N(x)\cap V(C)$.

Suppose that there is no feasible pair.
We say a vertex $x_j\in V(C)$ is of {\em type} $i$ if $N_H(x_j)\subseteq V(H_i-x)$ for some $i\in \{1,2\}$.
Then every vertex in $C$, except $x_0$ and $x_2$, must be of certain type.
By symmetry, let $x_{k-2}$ be of type 1, then we can infer (in order) that $x_{k-4},x_{k-6},...,x_1,x_{k-1},x_{k-3},...,x_4$ must be all of type 1,
and moreover $N_H(x_2)\subseteq V(H_1)$.
This shows that $\{x,x_0\}$ is a 2-cut of $G$ separating $H_2$ and $G-H_2$, but $G$ is 3-connected, a contradiction.

Hence there exist $x_i,x_{i+2}\in V(C)$ and $y\in V(H_1-x), z\in V(H_2-x)$ such that $x_iy,x_{i+2}z\in E(G)$.
Suppose that $N(x)\cap V(C)\neq \emptyset$.
By Claim 3, $x,y,z\in A$. So every $(y,z)$-path $P$ in $H$ passes through $x$ and thus is of even length at least 4.
Then $(C-\{x_{i+1}\})\cup x_iy \cup P\cup zx_{i+2}$ is an odd cycle of length at least $k+4$, a contradiction. This proves Claim 9.
\end{pf}

\begin{claim}
$|V(H)|=3$.
\end{claim}

\begin{pf}
By Claims 8 and 9, there exist a cut-vertex $x$ of $H$ with $N(x)\cap V(C)=\emptyset$
and a feasible pair $\{x_i,x_{i+2}\}$, where $V(H_1\cap H_2)=\{x\}$ and $H_1\cup H_2=H$.
Choose vertices $y_1\in N(x_i)\cap V(H_1-x)$ and $y_2\in N(x_{i+2})\cap V(H_2-x)$.
By Claim 2, $y_1,y_2\in A$.
If there is a $(y_1,y_2)$-path $P$ in $H$ with length at least $4$,
then $(C-\{x_{i+1}\})\cup x_iy_1\cup P\cup y_2x_{i+2}$
is an odd cycle of length at least $k+4$. So all $(y_1,y_2)$-paths in $H$ are of length 2.
This shows that for each $j\in \{1,2\}$, $y_jx\in E(G)$ and $H-y_jx$ is disconnected.
If $|V(H)|\ge 4$, then there is some $|V(H_j)|\ge 3$ and thus $y_j$ is a cut-vertex of $H$,
which is a contradiction to Claim 9. Thus $|V(H)|=3$.
\end{pf}

\medskip

By Claims 8 and 10, let $V(H)=\{x,z_1,z_2\}$ such that $xz_1,xz_2\in E(G)$ and $z_1z_2\notin E(G)$.
Claim 9 shows that $N_H(C)\subseteq \{z_1,z_2\}$.
So each vertex in $V(C)$ is adjacent to $z_1$ or $z_2$,
which will force triangles in $G$.
This contradiction completes the proof of Theorem \ref{Thm:kl}.{\hfill$\Box$}

\medskip
It remains to show the proofs of Lemmas \ref{Lem:3connected} and \ref{Lem:2oddcycle},
which we leave to Sections \ref{sec:3connected} and \ref{sec:2oddcycle}, respectively.

\section{Proof of Lemma \ref{Lem:3connected}}
\label{sec:3connected}
In this section, we establish Lemma \ref{Lem:3connected}, which we restate below for the
reader's convenience.

\medskip

\noindent {\bf Lemma \ref{Lem:3connected}.}
{\em
Let $G$ be
\begin{equation}
\label{equ:4critical}
\text{ a 4-critical graph with } L(G)=\{k,k+2l\}, \text{ where } k\geq 5 \text{ and } l\geq 1.
\end{equation}
Then $G$ is 3-connected.
}

\medskip

Clearly every graph $G$ satisfying \eqref{equ:4critical} is 2-connected with $\delta(G)\ge 3$.
The following weak version of Lemma \ref{Lem:2oddcycle} will be crucial in the proof of Lemma \ref{Lem:3connected}.

\begin{lemma}\label{Lem:2oddcycle-1}
For any graph $G$ satisfying \eqref{equ:4critical}, every two odd cycles intersect.
\end{lemma}

Let us first prove Lemma \ref{Lem:3connected}, assuming the above lemma.

\medskip

\noindent{\bf Proof of Lemma \ref{Lem:3connected}. (Assuming Lemma \ref{Lem:2oddcycle-1})}
The proof technique is similar to Corollary 4.2 in \cite{KRS}. Suppose that $G$ is not 3-connected.
Then there exists a 2-separator $(A,B)$ of $G$ such that $V(G)=A\cup B$,
$A\cap B=\{x,y\}$ and no edges are from $G[A]-\{x,y\}$ to $G[B]-\{x,y\}$.
We need a result from \cite{KRS} (see its Lemma 1.2), which states that
for any two vertices $v_1,v_2$ in a 4-critical graph, there is an odd cycle containing $v_1$ and avoiding $v_2$.
So for the vertex $x$ and any vertex $u\in A-\{x,y\}$, there is an odd cycle $C_1$ in $G$ containing $u$ and avoiding $x$;
and for the vertex $y$ and any vertex $v\in B-\{x,y\}$, there is an odd cycle $C_2$ in $G$
containing $v$ and avoiding $y$. It is easy to see that $V(C_1)\subseteq A-\{x\}$ and $V(C_2)\subseteq B-\{y\}$,
which imply that $V(C_1\cap C_2)=\emptyset$.
However $C_1$ and $C_2$ are odd, contradicting Lemma \ref{Lem:2oddcycle-1}.
The proof of Lemma \ref{Lem:3connected} is finished.{\hfill$\Box$}

\medskip

In the remainder of this section, we prove Lemma \ref{Lem:2oddcycle-1}.
To do so, as $L(G)=\{k,k+2l\}$, we consider three situations: (i) two $(k+2l)$-cycles; (ii) one $(k+2l)$-cycle and one $k$-cycle;
and (iii) two $k$-cycles. We will demonstrate each of the situations in a following separated subsection.

The next result will be used several times in this and forthcoming sections.

\begin{theorem}\emph{(\cite[Theorem~3.1]{KRS})}\label{THkrsextending}
Let $G$ be a graph with $|L(G)|=1$ and $C$ be an odd cycle in $G$. If
$G$ contains no $K_4$, then any proper 3-coloring of $C$ can be extended to a proper 3-coloring of $G$.
\end{theorem}

\subsection{$(k+2l)$-cycles intersect}
We first consider the case of two $(k+2l)$-cycles and show that it holds even for Lemma \ref{Lem:2oddcycle}.
\begin{lemma}
\label{Lem:Two(k+2l)}
For any graph $G$ satisfying \eqref{equ:4critical}, every two $(k+2l)$-cycles intersect in at least two vertices.
\end{lemma}

\noindent{\bf Proof.}
Suppose to the contrary that there exist two $(k+2l)$-cycles $C_0,C_1$ in $G$
with $|V(C_0\cap C_1)|\le 1$. Since $G$ is 2-connected, there are two disjoint $(C_0,C_1)$-paths,
say $R,S$, from $x_0,x_1\in V(C_0)$ to $y_0,y_1\in V(C_1)$, respectively.
In the case that $\lvert V(C_0 \cap C_1) \rvert=1$, we choose $R=V(C_0 \cap C_1)$.
So we always have $|S|\ge 1$.
Let $X$ be an $(x_0,x_1)$-path in $C_0$ and $Y$ a $(y_0,y_1)$-path in $C_1$
such that $C_2:=X\cup Y\cup R\cup S$ is an odd cycle.
Then $C_3:=(C_0 \cup C_1- \widehat{X}\cup \widehat{Y})\cup R\cup S$ is also an odd cycle.
But $|C_2|+|C_3|=2(|R|+|S|)+|C_0|+|C_1|>2(k+2l)$, a contradiction to $L(G)=\{k,k+2l\}$.
This proves the lemma.{\hfill$\Box$}

\subsection{$(k+2l)$-cycle intersects with $k$-cycle}
We then consider two odd cycles of different lengths.
\begin{lemma}
\label{Lem:(k+2l)+(k)}
For any graph $G$ satisfying \eqref{equ:4critical}, every $k$-cycle and $(k+2l)$-cycle intersect.
\end{lemma}

\noindent{\bf Proof.} Suppose to the contrary that there exist some $k$-cycle $C_0$ and $(k+2l)$-cycle $C_1$
in $G$ with $V(C_0\cap C_1)=\emptyset$. We will prove three claims, which lead us to contradictions.

\medskip

{\bf (A).}
For any vertex $u\in V(C_0)$, there is a $(u,C_1)$-path internally
disjoint from $C_0\cup C_1$.

\medskip

\begin{pf}
Since $C_0$ is induced (as it is a shortest odd cycle) and $\delta(G)\geq 3$,
for any vertex $u\in V(C_0)$, there exists a neighbor of $u$ not in $C_0$.
Now suppose that (A) fails. Then there exist some $u\in V(C_0)$ and
$C_0$-bridge $H$ such that $u\in V(H)$ and $V(C_1\cap H)=\emptyset$.
Let $G_0:=G[H\cup C_0]$ and $G_1:=G-V(H-C_0)$.
Note that $G_1$ is a proper subgraph of $G$.
Since $G$ is 4-critical, $G_1$ has a proper 3-coloring $\varphi$.
If there is a $(k+2l)$-cycle in $G_0$, say $C_2$, then
$V(C_1\cap C_2)\subset V(C_0\cap C_1)=\emptyset$,
a contradiction to Lemma \ref{Lem:Two(k+2l)}. Thus $L(G_0)=\{k\}$.
By Theorem \ref{THkrsextending}, the restriction of $\varphi$ on $C_0$ can be extended to
a proper 3-coloring of $G_0$. This gives a proper 3-coloring of $G$, a contradiction to \eqref{equ:4critical}.
\end{pf}

\medskip

{\bf (B).}
Let $R,S$ be any two disjoint $(C_0,C_1)$-paths from $x_0,x_1\in V(C_0)$
to $y_0,y_1\in V(C_1)$, respectively. Let $X$ be any path from $x_0$ to
$x_1$ on $C_0$, and $Y$ be any path from $y_0$ to $y_1$ on $C_1$. Then
$|R|+|S|=l$, and $|X|\in \{k+l-|Y|,|Y|-l\}$.

\medskip

\begin{pf}
Set $C_2:=X\cup Y\cup R\cup S$, and
$C_3:=(C_0\cup C_1-\widehat{X}\cup \widehat{Y})\cup R\cup S$.
If $C_2$ is odd, then $C_3$ is odd,
and $|C_2|+|C_3|=2(|R|+|S|)+2k+2l$.
Since $L(G)=\{k,k+2l\}$, we can then infer that $|C_2|=|C_3|=k+2l$,
$|R|+|S|=l$ and $|Y|=(k-|X|)+l$. If $C_2$
is even, repeat the above proof using $X':=C_0-X$ instead of $X$.
In this case, it holds that $|R|+|S|=l$ and $|Y|=(k-|X'|)+l$, implying that $|Y|=|X|+l$.
\end{pf}

\medskip

{\bf (C).}
There are three disjoint $(C_0,C_1)$-paths.

\medskip

\begin{pf}
Since $G$ is 2-connected, there are two disjoint
$(C_0,C_1)$-paths, say $R,S$, from $x_0,x_1\in V(C_0)$ to
$y_0,y_1\in V(C_1)$, respectively. By (B), $|R|+|S|=l$.
Let $P, Q$ be the two $(x_0,x_1)$-paths on $C_0$ with $|P|\le |Q|$.
Since $|C_0|\geq 5$, we have $|Q|\geq 3$.
Let $x_2$ be any vertex in $V(Q)\backslash\{x_0,x_1\}$.
We draw $C_0$ in the plane such that $x_0,x_1,x_2$ appear on $C$ clockwise.
Define $\alpha_i:=|x_iC_0x_{i+1}|$, where subscripts are taken modulo 3.
By (A), there is an $(x_2,C_1)$-path, say $T$, internally disjoint from $C_0\cup C_1$.
Suppose that (C) fails. Then every such $T$ intersects with $R\cup S$.

Let $z\in V(T\cap (R\cup S))$ such that $|x_2Tz|$ is minimal. If $z\in V(R)$,
then $R':=x_2Tz\cup zRy_0$ and $S$ are two disjoint
$(C_0,C_1)$-paths. Consider the following paths $x_0C_0x_1,x_1C_0x_2,x_2C_0x_1$ and $x_1C_0x_0$.
By (B), we obtain that $\alpha_0,\alpha_1,\alpha_0+\alpha_2,\alpha_1+\alpha_2\in \{k+l-|Y|,|Y|-l\}$, where $Y$ is a $(y_0,y_1)$-path
on $C_1$. Since $\alpha_2>0$, we have $\alpha_0=\alpha_1$.
If $z\in V(S)$, then by symmetry, we obtain $\alpha_0=\alpha_2$.
Note that $x_2$ can be picked to be any vertex in $V(Q)\backslash\{x_0,x_1\}$.
This shows that for any such $x_2$, either $|x_1C_0x_2|$ or $|x_2C_0x_0|$ equals $\alpha_0$.
Thus $|V(Q)\backslash\{x_0,x_1\}|\leq 2$, which, together with $|Q|\ge 3$, imply that $|Q|=3$.
Then $|C_0|=5$ and $|P|=2$.

Let $a$ be the vertex in $P-\{x_0,x_1\}$ and let $b,c$ be the vertices in $Q$ such that $Q=x_1bcx_0$. By (A), there exist a $(c,C_1)$-path $T_1$ and a $(b,C_1)$-path $T_2$, where each of them is internally disjoint from $C_0 \cup C_1$. Since we assume that (C) fails, each of $T_1,T_2$ intersects $R \cup S$. Applying the arguments in the previous paragraph with choosing $x_2=c$, we have that $T_1$ contains a subpath $cT_1u$ for some vertex $u$ in $R$ internally disjoint from $R \cup S$ since $\alpha_1=1$ and $\alpha_0=2$. Similarly, $T_2$ contains a subpath $bT_2v$ for some vertex $v$ in $S$ internally disjoint from $R \cup S$. Then $cT_1u \cup uRy_0$ and $bT_2v \cup vSy_1$ are disjoint $(C_0,C_1)$-paths. (Indeed, by a similar argument as above and (B), one can see $cT_1u$ and $bT_2v$ are disjoint.) Let $Y$ be a $(y_0,y_1)$-path on $C_1$. By (B),
$\{|x_0ax_1|,|x_1bcx_0|,|bc|\}=\{1,2,3\}\subseteq \{k+l-|Y|,|Y|-l\}$, which of course is a contradiction. This proves (C).
\end{pf}

\medskip

Hence, there are three disjoint $(C_0,C_1)$-paths $P_i$ from some vertex $x_i \in V(C_0)$ to some vertex $y_i \in V(C_1)$, for $i \in \{0,1,2\}$.
By (B), $|P_1|+|P_2|=|P_1|+|P_3|=|P_2|+|P_3|=l$, thus $3l=2(|P_1|+|P_2|+|P_3|)$, implying that
$l$ is even.

Observe that the subgraph $C_0\cup C_1\cup P_0\cup P_1\cup P_2$ is planar.
So we can draw it in the plane such that $x_0,x_1,x_2$ appear on $C_0$ clockwise and
$y_0,y_1,y_2$ appear on $C_1$ counterclockwise.
Define $\alpha_i:=|x_iC_0x_{i+1}|$ and $\beta_i:=|y_{i+1}C_1y_i|$, where the subscripts are taken
modulo 3. So $\alpha_0+\alpha_1+\alpha_2=k$ and
$\beta_0+\beta_1+\beta_2=k+2l$.
By (B), for any $i\in \{0,1,2\}$,
$\beta_i=\alpha_i+l$ or $\beta_i+\alpha_i=k+l$.
We discuss all possible cases.
If $\beta_i=\alpha_i+l$ for all $i$, then
$\beta_0+\beta_1+\beta_2=k+3l$, a contradiction.
If $\beta_i+\alpha_i=k+l$ for all $i$,
then $\beta_0+\beta_1+\beta_2=3k+3l-(\alpha_0+\alpha_1+\alpha_2)=2k+3l$,
a contradiction. If exactly two $i$'s satisfy $\beta_i=\alpha_i+l$, say $i=0,1$,
then $\beta_0+\beta_1+\beta_2=k+3l+\alpha_0+\alpha_1-\alpha_2$, implying
that $k+l=2\alpha_2$ is even, a contradiction to the facts that $k$ is odd and $l$ is even.
So there is exactly one $i$, say $i=0$, satisfying $\beta_i=\alpha_i+l$.
Then $\beta_0+\beta_1+\beta_2=2k+3l+\alpha_0-\alpha_1-\alpha_2$. This shows that $k+2l=k+3l+2\alpha_0$,
which cannot be. This finishes the proof of Lemma \ref{Lem:(k+2l)+(k)}.
{\hfill$\Box$}

\subsection{$k$-cycles intersect}
Lastly, we consider two $k$-cycles and prove Lemma \ref{Lem:Two(k)}, thereby completing the proof of Lemma \ref{Lem:3connected}.

Our proof is dependent of a well-known result due to Dirac \cite{D} (also see \cite[pp.367--368]{BM}).
Let $\{u,v\}$ be a 2-cut of a $k$-critical graph $G$ and $H$ be a component in $G-\{u,v\}$. We say $G[H\cup \{u,v\}]$ is of type 1
if every $(k-1)$-coloring of $G[H\cup \{u,v\}]$ assigns the same color to $u$ and $v$, and of type 2 if every $(k-1)$-coloring of $G[H\cup \{u,v\}]$ assigns distinct colors to $u$ and $v$.

\begin{theorem}\cite{D}
\label{Thm:Dirac}
Let $G$ be a $k$-critical graph with a 2-vertex cut $\{u,v\}$. Then there exists a 2-separation of $G$, say $(V_1,V_2)$,
with $V_1\cap V_2=\{u,v\}$, such that:\\
(1) $uv\notin E(G)$;\\
(2) $G=G_1\cup G_2$, where $G_i=G[V_i]$ is of type $i$, for $i=1,2$;\\
(3) both $G_1+uv$ and $G_2/\{u,v\}$ are $k$-critical.
\end{theorem}

\begin{lemma}
\label{Lem:Two(k)}
For any graph $G$ satisfying \eqref{equ:4critical}, every two $k$-cycles intersect.
\end{lemma}

\noindent{\bf Proof.}
Suppose to the contrary that there exist two $k$-cycles $C_0,C_1$ in $G$ with $V(C_0\cap C_1)=\emptyset$.
The case $l=1$ was solved in Proposition 4.1 of \cite{KRS}, so we assume that $l\ge 2$.
Write $C_0:=x_0x_1...x_{k-1}x_0$ throughout this proof.
We divide the proof into a sequence of claims.

\setcounter{claim}{0}
\begin{claim}
For each $i\in \{0,1\}$ and each vertex $u\in V(C_i)$, there is a $(u,C_{1-i})$-path $P_u$ in
$G$, internally disjoint from $C_0\cup C_1$.
\end{claim}

\begin{pf}
By symmetry, we may only consider vertices in $C_0$.
Suppose to the contrary that there exists $u\in C_0$ and some $C_0$-bridge $H$
such that $u\in V(H)$ and $V(H\cap C_1)=\emptyset$.
Let $G_0:=G[H\cup C_0]$ and $G_1:=G-V(H-C_0)$.
As $G$ is 4-critical, $G_1$ has a proper 3-coloring $\varphi$. If $G_0$ contains a
$(k+2l)$-cycle, say $C_2$, then
$V(C_1\cap C_2)\subset V(C_0\cap C_1)=\emptyset$,
a contradiction to Lemma \ref{Lem:(k+2l)+(k)}. Thus $L(G_0)=\{k\}$.
By Theorem \ref{THkrsextending}, the restriction of $\varphi$ on $C_0$ can be extended to a proper $3$-coloring of $G_0$.
This gives a proper 3-coloring of $G$, a contradiction.
\end{pf}

\begin{claim}
Let $R,S$ be any two disjoint $(C_0,C_1)$-paths from
$x_i,x_j\in V(C_0)$ to $y_i,y_j\in V(C_1)$, respectively.
Let $X$ be any $(x_i,x_j)$-path on $C_0$ and $Y$
be any $(y_i,y_j)$-path on $C_1$. Let $t:=|R|+|S|$.
Then $t\in \{l,2l\}$. If $t=l$, then $||X|+|Y|-k|=l$ or $||X|-|Y||=l$; and if $t=2l$,
then $|Y|=|X|$ or $k-|X|$.
In particular, when $1\le |X|\le 2$, we have $|Y|\in \{l+|X|,k-l-|X|\}$ if
$t=l$, and $|Y|\in \{|X|,k-|X|\}$ if $t=2l$.
\end{claim}

\begin{pf}
Let $C_2:=X\cup Y\cup R\cup S$ and $C_3:=(C_0\cup C_1-\widehat{X}\cup \widehat{Y})\cup R\cup S$.
Then $|C_2|+|C_3|=2(|R|+|S|)+|C_0|+|C_1|=2(|R|+|S|)+2k$.
If $C_2$ is odd, then $C_3$ is also odd. There are three cases:
($a$) $|C_2|=k$, $|C_3|=k+2l$; ($b$) $|C_2|=k+2l$, $|C_3|=k$;
and ($c$) $|C_2|=|C_3|=k+2l$. In each case, we can infer that
$|Y|=(k-|X|)-l$ and $t=l$; $|Y|=(k-|X|)+l$ and $t=l$;
$|Y|=k-|X|$ and $t=2l$, respectively.
Otherwise, $C_2$ is even. Then we can repeat the above proof by using $X':=C_0-X$ instead of $X$.
Similarly, we have $t\in \{l,2l\}$, and it is a routine matter to verify other quantities.
The result when $1\le |X|\le 2$ easily follows by the facts that $l\ge 2$ and $1\le |Y|\le k-1$.
\end{pf}

\begin{claim}
Let $P_1,P_2,P_3$ be three disjoint $(C_0,C_1)$-paths of $G$, with $|P_1|\leq |P_2|\leq |P_3|$.
Then one of the followings holds:\\
(1) $|P_1|=|P_2|=|P_3|=l/2$;
(2) $|P_1|=|P_2|=|P_3|=l$;
(3) $|P_1|=|P_2|=l/2,|P_3|=3l/2$.
\end{claim}

\begin{pf}
Set $|P_1|=a,|P_2|=b,|P_3|=c$.
By Claim 2, each of $a+b,a+c$ and $b+c$ must be in $\{l,2l\}$.
Consider the vector $(a+b,a+c,b+c)$, which cannot be $(l,l,2l)$.
Therefore, the vector only can be $(l,l,l), (l,2l,2l)$ or $(2l,2l,2l)$, which gives that
$(a,b,c)=(\frac{l}{2},\frac{l}{2},\frac{l}{2}), (\frac{l}{2},\frac{l}{2},\frac{3l}{2})$ or $(l,l,l)$,
respectively. This proves Claim 3.
\end{pf}

\begin{claim}
There exist two disjoint $(C_0,C_1)$-paths from two consecutive vertices $x_i,x_{i+1}\in V(C_0)$ to $V(C_1)$
for some $i \in \{0,1,...,k-1\}$.
\end{claim}

\begin{pf}
Since $G$ is 2-connected, there are two disjoint $(C_0,C_1)$-paths $P_1,P_2$,
say from $x_i,x_j\in V(C_0)$ to $y_i,y_j\in V(C_1)$, respectively.
We choose $P_1,P_2$ such that $d_{C_0}(x_i,x_j)$ is minimal.
It is enough to show that $d_{C_0}(x_i,x_j)=1$.
Suppose to the contrary that there exists some vertex $x_m\in V(X)-\{x_i,x_j\}$,
where $X$ is the shorter $(x_i,x_j)$-path on $C_0$.
By Claim 1, there exists an $(x_m,C_1)$-path $Q$, which is internally disjoint from $C_0\cup C_1$.
If $Q$ is disjoint from some $P_t$, then $P_t,Q$ is a pair of
disjoint $(C_0,C_1)$-paths with a shorter distance on $C_0$,
a contradiction. So $Q$ intersects $P_1\cup P_2$.
Let $z\in V(Q)\cap V(P_1\cup P_2)$ be the vertex such that $|x_mQz|$ is the minimum,
say $z\in P_1$. Then $x_mQz\cup zP_1y_1$ together with $P_2$
form a pair of disjoint $(C_0,C_1)$-paths such that the length of the
shortest path in $C_0$ connecting their ends in $C_0$ is less than $d_{C_0}(x_i,x_j)$, a contradiction.
\end{pf}

\medskip

Let $P_i,P_{i+1}$ be two disjoint $(C_0,C_1)$-paths from consecutive $x_i,x_{i+1}\in V(C_0)$
to some $y_i,y_{i+1}\in V(C_1)$, respectively.
If exist, let $P_{i+2}$ be a path from $x_{i+2}$ to $z\in V(P_i\cup P_{i+1})-\{x_i,x_{i+1}\}$ internally
disjoint from $P_i\cup P_{i+1}\cup C_0\cup C_1$.
If such $P_{i+2}$ does not exist, then the coming Claim 6 will hold trivially and in this case readers can skip Claim 5 and the proof of Claim 6.
Let $t:=|P_i|+|P_{i+1}|$.

\begin{claim} Assume that $P_{i+2}$ exists. (1) If $z\in V(P_{i+1})$, then $k=l+3$, $|x_{i+1}P_{i+1}z|=t-l+1$ and $|P_{i+2}|=2l+1-t$.\\
(2) If $z\in V(P_i)$, then $|x_iP_iz|=|P_{i+2}|=1$ or $l+1$; in the latter case, we have $z=y_i$.
\end{claim}

\begin{pf}
Without loss of generality, let $i=0$ and $P_0,P_1$ be from $x_0,x_1$ to $y_0,y_1\in V(C_1)$, respectively.
By Claim 2, $t\in \{l,2l\}$.
Let $Y$ be a $(y_0,y_1)$-paths of $C_1$.

First consider $z\in V(P_1)$. Let $P'_2:=P_2\cup zP_1y_1$ and $C_2:=x_1P_1z\cup P_2\cup x_1x_2$.
Let $s:=|P_0|+|P'_2|$.
By Claim 2, if $s=l$ then $|Y|\in \{l+2,k-l-2\}$; and if $s=2l$ then
$|Y|\in \{2,k-2\}$. Similarly, if $t=l$ then $|Y|\in \{l+1,k-l-1\}$; and if $t=2l$ then
$|Y|\in \{1,k-1\}$.
As a consequence, $t,s$ cannot both be $2l$ (as, otherwise, we can obtain $k=3$, a contradiction).
If $t=s=l$, then $\{l+1,k-l-1\}\cap \{l+2,k-l-2\}\neq \emptyset$ implies that
$k=2l+3$; moreover, $x_1P_1z$ has the same length as $P_2$, implying that $C_2$ is odd.
Thus $|x_1P_1z|+|P_2|\geq k-1=2l+2$, contradicting the fact that $|x_1P_1z|+|P_2|\leq |P_1|+|P'_2|\leq s+t=2l$.
Hence, $\{t,s\}=\{l,2l\}$, and in this case, we can always get $k=l+3$.
Let $r:=\min \{|x_1P_1z|,|P_2|\}$ and $r':=\max\{|x_1P_1z|,|P_2|\}$.
Note that $|C_2|=1+(r'-r)+2r=1+|t-s|+2r=k-2+2r$ is odd.
If $|C_2|=k+2l$, then $r=l+1\le \min \{t,s\}=l$, a contradiction. So $|C_2|=k$ and thus $r=1$, $r'=1+l$.
The left part is easy to check. This proves (1).

Now suppose $z\in V(P_0)$.
Let $P'_2:=P_2\cup zP_0y_0$ and $s:=|P_1|+|P'_2|$.
By Claim 2, $s,t\in \{l,2l\}$. If $s\neq t$,
then by Claim 2, $|Y|\in \{1,k-1\}\cap \{l+1,k-l-1\}$.
This implies $k=l+2$. Let $r:=\max \{|x_0P_0z|,|P_2|\}$
and $r':=\min \{|x_0P_0z|,|P_2|\}$ with $r-r'=|s-t|=l$.
Then $x_0x_1x_2\cup P_2\cup x_0P_0z$ is an odd cycle
with length $2+l+2r'=k+2r'$, implying that $r'=l$ and $r=2l$. This is a contradiction
to $r<\max\{s,t\}=2l$. So $s=t$. Then $|x_0P_0z|=|P_2|$, implying that
$C_3:=x_0P_0z\cup P_2\cup (C_0-\{x_1\})$ is an odd cycle with length $k-2+2|P_2|$.
Thus $|P_2|= 1$ or $l+1$. In the later case, $C_3$ is of length $k+2l$
and by Lemma \ref{Lem:(k+2l)+(k)}, we must have $z=y_0$. This proves (2).
\end{pf}

\begin{claim}
There exist three disjoint $(C_0,C_1)$-paths from consecutive $x_i,x_{i+1},x_{i+2}\in V(C_0)$ to $V(C_1)$.
\end{claim}

\begin{pf}
By Claim 4, we may assume that there exist two disjoint $(C_0,C_1)$-paths $P_0,P_1$
from $x_0,x_1$ to $y_0,y_1\in V(C_1)$, respectively. Write $t:=|P_0|+|P_1|\in \{l,2l\}$ (by Claim 2).
For each $i\in \{2,k-1\}$, let $P_i'$ be a $(x_i,C_1)$-path satisfying the conclusion of Claim 1.
Suppose that each $P_i'$ intersects $P_0\cup P_1$. Let $z$ be the vertex in $V(P_2') \cap V(P_0 \cup P_1)$
such that the path $x_2P_2'z$, denoted by $P_2$, is as short as possible. Similarly, let $w$ be the
vertex in $V(P_{k-1}') \cap V(P_0 \cup P_1)$ such that the path $x_{k-1}P_{k-1}'w$, denoted by
$P_{k-1}$, is as short as possible.

We show that $P_2$ and $P_{k-1}$ are internally disjoint. Suppose not.
Then there exists $x\in V(P_2\cap P_{k-1})$ such that $xP_{k-1}x_{k-1}, xP_2x_2$ and
$xP_2z$ are internally disjoint. At this point, one actually need not to distinguish between $xP_2z$ and $xP_{k-1}w$,
and thus we may assume, without loss of generality, that $z\in V(P_1)$.
Let $C_2:=x_2P_2x\cup xP_{k-1}x_{k-1}\cup (C_0-\{x_0,x_1\})$,
and $C_3:=x_2P_2x\cup xP_{k-1}x_{k-1}\cup (x_{k-1}x_0x_1x_2)$.
If $C_2$ is odd, as $V(C_2\cap C_1)=\emptyset$, by Lemma \ref{Lem:(k+2l)+(k)} we infer that $|C_2|=k$.
Let $P_3:=xP_2z\cup zP_1y_1$, $P_4:=(x_2x_1x_0)\cup P_0$, and $P_5:=x_{k-1}x_0\cup P_0$.
Note that $\{P_3,P_4\}$ and $\{P_3,P_5\}$ are both pairs of disjoint $(C_1,C_2)$-paths.
By Claim 2, $|P_3|+|P_4|, |P_3|+|P_5|\in \{l,2l\}$, where $|P_5|=|P_4|+1$.
This implies $l=1$, a contradiction.
Hence $C_3$ is odd.
As $V(C_3\cap C_2)=\emptyset$, we infer that $|C_3|=k$.
Note that $P_0,P_3$ are disjoint $(C_1,C_3)$-paths, and $P_0, P_2\cup zP_1y_1$
are disjoint $(C_0,C_1)$-paths. Since $P_3$ is a subpath of $P_2\cup zP_1y_1$,
we have $|P_0|+|P_3|=l$ and $|P_0|+|P_2\cup zP_1y_1|=2l$, implying $|x_2P_2x|=l$.
By Claim 5, $k=l+3$. Then $|C_3|>l+3=k$,
a contradiction. Therefore indeed $P_2$ and $P_{k-1}$ are internally disjoint.

Next we discuss the locations of $z,w$. First assume that both $z,w\in V(P_j)$, say $j=1$.
By Claim 5(1), we have $k=l+3$, $|x_1P_1z|=t-l+1$ and $|P_2|=2l+1-t$;
and by Claim 5(2), we get
$|x_1P_1w|=|P_{k-1}|=1$ or $l+1$.
Let $C_4:=P_2\cup P_{k-1}\cup wP_1z\cup (x_2x_1x_0x_{k-1})$.
If $t=l$ and $|x_1P_1w|=|P_{k-1}|=1$, then $|x_1P_1z|=1$ (implying $z=w$) and $|P_2|=l+1=k-2$,
implying $|C_4|=k+2$ and thus $l=1$, a contradiction.
If $t=l$ and $|x_1P_1w|=|P_{k-1}|=l+1$,
then $|x_1P_1z|=1$, $|P_2|=l+1$ and $w=y_1$ (by Claim 5(2)),
implying that $|C_4|=2(l+1)+l+3=k+2l+2$, a contradiction. Hence $t=2l$.
So $|x_1P_1z|=l+1$ and $|P_2|=1$.
If $|x_1P_1w|=|P_{k-1}|=1$,
then $w\in x_1P_1z$ and $|wP_1z|=l$, implying $|C_4|=l+5=k+2$, a contradiction;
otherwise $|x_1P_1w|=|P_{k-1}|=l+1$, then $w=z$, also implying $|C_4|=l+5=k+2$, a contradiction.

Suppose $z\in V(P_1)$ and $w\in V(P_0)$.
By Claim 5, $k=l+3$, $|x_1P_1z|=|x_0P_0w|=t-l+1$ and $|P_2|=|P_{k-1}|=2l+1-t$.
Then $(C_0-\{x_0x_{k-1},x_1x_2\})\cup x_1P_1z\cup P_2\cup x_0P_0w\cup P_{k-1}$
is an odd cycle of length $(k-2)+2(l+2)=k+2l+2$, a contradiction.

Lastly we consider $z\in V(P_0)$ and $w\in V(P_1)$. By Claim 5,
$|P_2|=|x_0P_0z|\in \{1,l+1\}$ and $|P_{k-1}|=|x_1P_1w|\in \{1,l+1\}$.
Then $(C_0-\{x_0x_{k-1},x_1x_2\})\cup x_1P_1w\cup P_{k-1}\cup x_0P_0z\cup P_2$
is an odd cycle of some length $s$, where $s\in \{k+2,k+2l+2,k+4l+2\}$. Note that
each case yields a contradiction.
This completes the proof of Claim 6.
\end{pf}

\medskip

In the rest of this proof, we write $C_1=u_0u_1...u_{k-1}u_0$.
By Claim 6, we may assume that there exist three disjoint $(C_0,C_1)$-paths $P_0,P_1,P_2$ from
consecutive $x_0,x_1,x_2\in V(C_0)$ to $y_0,y_1,y_2\in V(C_1)$, respectively.
In view of Claim 3, we can get
\begin{align*}
(|P_0|,|P_1|,|P_2|)\in \left\{(\frac{l}{2},\frac{3l}{2},\frac{l}{2}),(\frac{3l}{2},\frac{l}{2},\frac{l}{2}),
(\frac{l}{2},\frac{l}{2},\frac{3l}{2}),(l,l,l),(\frac{l}{2},\frac{l}{2},\frac{l}{2})\right\}.
\end{align*}
For each $i \in \{0,1,2\}$, define $Y_i$ to be the path of $C_1$ from $y_i$ to $y_{i+1}$ not
containing $y_{i+2}$, where the indices are taken modulo 3. Let $\beta_i:=|Y_i|$.
So $\beta_0+\beta_1+\beta_2=k$. Without loss of generality, we draw $C_0,C_1$
on the plane such that $x_1,x_2,x_3$ appear in $C_0$ in the clockwise
order and $u_0,u_1,...,u_{k-1}$ appear in $C_1$ in the counterclockwise order.

\begin{claim}
$(|P_0|,|P_1|,|P_2|)=(\frac{l}{2},\frac{l}{2},\frac{l}{2})$ and thus $l$ is even.
\end{claim}
\begin{pf}
First suppose $(|P_0|,|P_1|,|P_2|)=(\frac{l}{2},\frac{3l}{2},\frac{l}{2})$.
Note that $l$ is even. By Claim 2, $\beta_0,\beta_1\in \{1,k-1\}$ and $\beta_2\in \{l+2,k-l-2\}$.
As $\beta_1+\beta_2+\beta_3=k$, we must have $\beta_0=\beta_1=1$ and thus $\beta_2=k-2$ is odd,
a contradiction.

Next, we assume that $(|P_0|,|P_1|,|P_2|)=(\frac{3l}{2},\frac{l}{2},\frac{l}{2})$.
Note that $l$ is even.
By Claim 2, $\beta_0\in \{1,k-1\}$, $\beta_1\in \{l+1,k-l-1\}$ and $\beta_2\in \{2,k-2\}$.
Note that $\beta_0+\beta_1+\beta_2=k$. Clearly $\beta_0=1$.
If $\beta_2=k-2$, then $\beta_1$ is odd and thus $\beta_1=l+1$,
implying $\beta_0+\beta_1+\beta_2=k+l$, a contradiction.
Therefore $\beta_0=1$ and $\beta_2=2$. So $\beta_1$ is even and thus $\beta_1=k-l-1$.
This shows that $k=k-l+2$, implying $l=2$.
Now we have $|P_0|=3, |P_1|=|P_2|=|Y_0|=1$, $|Y_1|=k-3$ and $|Y_2|=2$.
Without loss of generality, let $y_0=u_0, y_1=u_1$ and $y_2=u_{k-2}$.
By Claim 1, there exists an $(x_3,C_1)$-path.
So there exists a path $P_3$ from $x_3$ to $z\in V(P_0\cup P_1\cup P_2\cup C_1)$
internally disjoint from $P_0\cup P_1\cup P_2\cup C_0\cup C_1$.
We consider the location of $z$.
If $z\in V(P_0)$, then $P_3\cup zP_0u_0$ and $P_2$ are two
disjoint $(C_0,C_1)$-paths from $x_3,x_2$ to $u_0,u_{k-2}$, respectively.
By Claim 2, $2=|Y_2|\in \{1,k-1\}$ or $\{l+1,k-l-1\}$. Since $l=2$, the only possibility is $2=k-3$.
Thus $k=5$ and $|P_3\cup zP_0u_0|+|P_2|=l=2$, implying
that $z=u_0$ and $|P_3|=1$. Then $P_3\cup P_0\cup (x_0x_1x_2x_3)$
is an odd cycle of length $7$, a contradiction.
So, $z\in V(C_1)-\{u_0\}$ (as $|P_1|=|P_2|=1$) and $P_3$ is disjoint from $P_0$.
By Claim 2, we see that $|P_0|+|P_3|\in \{2, 4\}$ and so $|P_3|=1$.
If $z=u_1$, then $P_3\cup (C_1-\{u_0u_1\})\cup P_0\cup (x_0x_1x_2x_3)$
is an odd cycle of length $k+6$, a contradiction.
If $z=u_{k-2}$, then $G$ has a triangle on $\{x_2,x_3,z\}$, a contradiction.
If $z=u_{k-1}$, then $P_3\cup (C_1-\{u_{k-1}u_{k-2}\})\cup P_2\cup x_2x_3$
is an odd cycle of length $k+2$, a contradiction to $l=2$.
Thus $z\in V(Y_1)-\{u_1,u_{k-2}\}$, then by Claim 2, $|zY_1u_{k-2}|\in \{3,k-3\}$.
Since $|zY_1u_{k-2}|<|Y_1|=k-3$, we obtain that $|zY_1u_{k-2}|=3$.
Then $P_3\cup zY_1u_1\cup P_1\cup (x_1x_2x_3)$ is an odd cycle of length $(k-6)+4=k-2$,
a contradiction. By symmetry, we can prove $(|P_0|,|P_1|,|P_2|)\neq(\frac{l}{2},\frac{l}{2},\frac{3l}{2})$.

Lastly, we suppose $(|P_0|,|P_1|,|P_2|)=(l,l,l)$. So $|P_i|+|P_j|=2l$.
By Claim 2, $\beta_0, \beta_1\in \{1,k-1\}$ and $\beta_2\in \{2,k-2\}$.
It is easy to see that $\beta_0=\beta_1=1$ and $\beta_2=k-2$.
Without loss of generality, let $y_i=u_i$ for $i\in \{0,1,2\}$.
By Claim 1, there exists an $(x_3,C_1)$-path internally disjoint from $C_0\cup C_1$.
So there exists a path $P_3$ from $x_3$ to
$z\in V(P_0\cup P_1\cup P_2\cup C_1)$ internally
disjoint from $P_0\cup P_1\cup P_2\cup C_0\cup C_1$.
Similarly as the above analysis, we consider four cases.
If $z\in V(P_0)$, then $P_3\cup zP_0u_0$ and $P_1$
are two disjoint $(C_0,C_1)$-paths from $x_3,x_1$ to $u_0,u_1$, respectively.
However, this is a contradiction to Claim 2,
as $|P_3\cup zP_0u_0|+|P_1|$ is larger than $l$
and thus equals $2l$, which implies $\beta_0\in \{2,k-2\}$.
If $z\in V(P_1)$, by Claim 5(2), $|x_1P_1z|=|P_3|$.
Then $(x_0x_1x_2x_3)\cup P_3\cup zP_1y_1\cup (C_1-Y_0)\cup P_0$
is an odd cycle of length $3+2l+(k-1)=k+2l+2$, a contradiction.
If $z\in V(P_2)$, by Claim 5(1), we get $k=l+3$, $|x_2P_2z|=l+1>|P_2|=l$, a contradiction.
Lastly, we consider $z\in V(C_1)-\{u_0,u_1,u_2\}$.
In this case, $P_1,P_2,P_3,P_4$ are four disjoint $(C_0,C_1)$-paths.
By Claim 3, $|P_3|=l$. By Claim 2, we see that $z=y_3$.
Then $(C_1-\{y_0y_1,y_2y_3\})\cup \{x_0x_1,x_2x_3\}\cup P_0\cup P_1\cup P_2\cup P_3$
is an odd cycle of length $k+4l$, a contradiction.
This proves Claim 7.
\end{pf}

\begin{claim}
$L(G)=\{5,9\}$, and any two disjoint 5-cycles $C_0,C_1$
in $G$ induce a Petersen graph $G[V(C_0\cup C_1)]$.
\end{claim}

\begin{pf}
Note that $l$ is even.
By Claim 2, $\beta_0,\beta_1\in \{l+1,k-l-1\}$ and $\beta_2\in \{l+2,k-l-2\}$.
Since $\beta_0+\beta_1\neq k$, we have $\beta_0=\beta_1$
and thus $\beta_2$ must be odd, so $\beta_2=k-l-2$.
Since $2(l+1)+(k-l-2)>k$, we have $\beta_0=\beta_1=k-l-1$ and thus $2(k-l-1) +(k-l-2)=k$.
We then get
\begin{equation}
\label{equ:k=3l/2+2}
k=\frac{3l}{2}+2,
\end{equation}
which implies that
$\beta_0=\beta_1=\frac{l}{2}+1$ and $\beta_2=\frac{l}{2}$.
Observe that $\frac{l}{2}$ is odd. Applying Claim 1 for $x_3$, we see that there is a path, say $P_3$,
from $x_3$ to $z\in V(P_0\cup P_1\cup P_2\cup C_1)$, internally disjoint from
$P_0\cup P_1\cup P_2\cup C_0\cup C_1$.

We show that $z\in V(C_1)-\{y_0,y_1,y_2\}$. Suppose not, we consider three cases that $z\in V(P_i)$ for $i=0,1,2$.
If $z\in V(P_0)$, then $P_3\cup zP_0y_0$ and $P_2$ are two disjoint $(C_0,C_1)$-paths from
$x_3,x_2$ to $y_0,y_2$, respectively. Let $t=|P_3\cup zP_0u_0|+|P_2|$. By Claim 2, if $t=l$,
then $\beta_0+\beta_1\in \{l+1,k-l-1\}=\{l+1,\frac{l}{2}+1\}$. However, $\beta_0+\beta_1=l+2$, a contradiction.
If $t=2l$, then $\beta_0+\beta_1\in \{1,k-1\}=\{1,\frac{3l}{2}+1\}$,
and thus $\frac{3l}{2}+1=l+2$, implying $l=2$. So, $k=5$ and $L(G)=\{5,9\}$. And $P_0:=x_0y_0\in E(G)$, $z=y_0$.
But $C_2:=P_3\cup P_0\cup (x_0x_1x_2x_3)$ gives a 7-cycle, yielding a contradiction. If $z\in V(P_1)$,
by Claim 5(2), $|P_3|=|x_1P_1z|$, and then, $P_3\cup zP_1y_1\cup Y_0\cup P_0\cup (x_0x_1x_2x_3)$
is an odd cycle of length $\frac{l}{2}+\frac{l}{2}+1+\frac{l}{2}+3=k+2$. This implies that $l=1$,
a contradiction. If $z\in V(P_2)$, then by Claim 5(1), we get $k=l+3$.
Together with \eqref{equ:k=3l/2+2}, this yields that $(l,k)=(2,5)$. By Claim 2, we obtain $|P_2|=1$, $|P_3|=3$ and
$z=y_2$. However, $P_3\cup zC_1y_0\cup P_0\cup (x_0x_1x_2x_3)$ is an odd cycle of length 11,
a contradiction. Therefore, indeed, $z\in V(C_1)-\{y_0,y_1,y_2\}$ and thus $P_0,P_1,P_2,P_3$ are pairwise disjoint paths.

By Claim 3, $|P_3|\in \{\frac{l}{2},\frac{3l}{2}\}$.
Let $y_0=u_0$, $y_1=u_{1+l/2}$, $y_2=u_{l+2}$, and $Y$ be a $(y_2,z)$-path on $C_1$.
First suppose that $|P_3|=\frac{3l}{2}$.
By Claim 2, we deduce that $|Y|\in \{1,k-1\}$, which implies $z=u_{l+1}$ or $u_{l+3}$. If $z=u_{l+3}$, then $P_3\cup y_0C_1z\cup P_0\cup (x_0x_1x_2x_3)$
is an odd cycle of length $\frac{3l}{2}+\frac{l}{2}-1+\frac{l}{2}+3=k+l$, a contradiction; if $z=u_{l+1}$,
then $P_3\cup zC_1y_1\cup P_1\cup (x_1x_2x_3)$ is also an odd cycle of length $\frac{3l}{2}+\frac{l}{2}+\frac{l}{2}+2=k+l$,
again a contradiction. Thus, $|P_3|=\frac{l}{2}$. From Claim 2, we infer that $z=y_1$ (this cannot happen)
or $z=u_1$. Then $(C_0-\{x_1,x_2\})\cup P_3\cup y_0z\cup P_0$
is an odd cycle of length $k+l-2$. This shows that $l=2$ and by \eqref{equ:k=3l/2+2}, we have $L(G)=\{5,9\}$.
We then see that $C_0,C_1$ are both 5-cycles, and
$x_0u_0,x_1u_2,x_2u_4,x_3u_1\in E(G)$. Consider the path $P_{4}$
from $x_{4}$ to $w\in V(P_0\cup P_1\cup P_2\cup C_1)$, internally disjoint from
$P_0\cup P_1\cup P_2\cup C_0\cup C_1$. By the symmetry between $x_3$ and $x_{4}$, similarly as above,
we can derive that $|P_{4}|=1$ and $w=u_3$.
In view of $L(G)=\{5,9\}$, now it is easy to verify that $G[V(C_0\cup C_1)]$ induces a Petersen graph.
\end{pf}

\begin{claim}
Let $H:=G[V(C_0\cup C_1)]$. For any $u\notin V(H)$, there are 3 disjoint paths from $u$ to $H$.
\end{claim}
\begin{pf}
Otherwise, there exists a 2-separation $(G_1,G_2)$ such that $G=G_1\cup G_2$,
$V(G_1\cap G_2)=\{a,b\}$, $\{u\}\subset G_1$, and $H\subset G_2$. Since $G$ is 4-critical,
by Theorem \ref{Thm:Dirac}, we have $ab\notin E(G)$ and either $G_1+ab$ or $G_1/\{a,b\}$ is 4-critical.
First we claim that there is a 5-cycle $D$ in $G_2$ which is disjoint with $\{a,b\}$.
This can be deduced from an easy observation that
Petersen graph has a 5-cycle disjoint from any two prescribed nonadjacent vertices of it.
Next we show that $G_1-a$ contains a 5-cycle. Note that in either case that $G_1+ab$ or $G_1/\{a,b\}$ is 4-critical,
we have that $G_1-a$ is 3-chromatic.
So, there always is an odd cycle in $G_1-a$, say $D'$. If $D'$ is a 9-cycle, then $D$ and $D'$ are
disjoint 5-cycle and 9-cycle in $G$, a contradiction to Lemma \ref{Lem:(k+2l)+(k)}. Thus $D'$ is a 5-cycle.
Observe that there are no two disjoint edges connecting $D$ and $D'$. So $G[V(D\cup D')]$ cannot be
a Petersen graph, a contradiction to Claim 8.
This proves Claim 9.
\end{pf}

\medskip

Now we are ready to finish the proof of Lemma \ref{Lem:Two(k)}. If $G=H$, then $G$ is 3-colorable,
a contradiction. Thus, there exists $u\notin V(H)$. By Claim 9, there are 3 disjoint paths, say
$Q_1,Q_2,Q_3$, from $u$ to $v_1,v_2,v_3\in V(H)$, respectively.
Since $G$ contains no triangle, we may assume that $v_1v_2\notin E(G)$.
Observe that there exist $(v_1,v_2)$-paths of lengths
2,3,4,5 in the Petersen graph $H$ for any nonadjacent vertices $v_1,v_2$.
These paths, together with the path $P_1\cup P_2$ (internally disjoint from $H$),
form two odd cycles of lengths differ by two, a contradiction to $L(G)=\{5,9\}$.
This final contradiction proves Lemma \ref{Lem:Two(k)}.{\hfill$\Box$}
\medskip

Putting Lemmas \ref{Lem:Two(k+2l)},~\ref{Lem:(k+2l)+(k)} and \ref{Lem:Two(k)} together,
we now complete the proof of Lemma \ref{Lem:2oddcycle-1}.

\section{Proof of Lemma \ref{Lem:2oddcycle}}
\label{sec:2oddcycle}
We devote this section to the proof of Lemma \ref{Lem:2oddcycle}, which asserts that
for any graph $G$ satisfying \eqref{equ:4critical}, every two odd cycles in $G$ intersect in at least two vertices.

In view of Lemma \ref{Lem:Two(k+2l)}, it is enough to consider two cases: (i) one $(k+2l)$-cycle and one $k$-cycle;
and (ii) two $k$-cycles.
To this end, we need the following two lemmas.

\begin{lemma}
\label{Lem:(k+2l)+(k)2}
For any graph $G$ satisfying \eqref{equ:4critical}, every $k$-cycle and $(k+2l)$-cycle intersect in at least two vertices.
\end{lemma}

\begin{lemma}
\label{Lem:Two(k)2}
For any graph $G$ satisfying \eqref{equ:4critical}, every two $k$-cycles intersect in at least two vertices.
\end{lemma}

\subsection{Proof of Lemma \ref{Lem:(k+2l)+(k)2}.}
Let $C_0$ be a $k$-cycle and $C_1$ be a $(k+2l)$-cycle. Suppose that
$C_0,C_1$ intersect in at most one vertex. By Lemma \ref{Lem:2oddcycle-1}, we
see $V(C_0\cap C_1)\neq \emptyset$. In the following, we denote the unique vertex in $C_0\cap C_1$ by $o$.

\setcounter{claim}{0}
\begin{claim}
For any $u\in C_0-\{o\}$, there is a path $P_u$ from $u$ to $u'\in C_1-\{o\}$, internally
disjoint from $C_0\cup C_1$.
\end{claim}
\begin{pf}
Suppose to the contrary that there exist $u\in C_0-\{o\}$ and some $C_0$-bridge
$H$ such that $u\in V(H)$ and $V(H\cap C_1)\subseteq \{o\}$.
Let $G_0:=G[H\cup C_0]$ and $G_1:=G-(H-C_0)$. Then $G_1$ has
a proper 3-coloring $\varphi$.
By Lemma \ref{Lem:Two(k+2l)}, we have $L(G_0)=\{k\}$.
Then Theorem \ref{THkrsextending} shows that the restriction $\varphi$ on $C_0$ can be extended
to a proper 3-coloring of $G_0$. Now this gives rise to a proper 3-coloring of $G$,
a contradiction.
\end{pf}

\begin{claim}
Let $P$ be a $(C_0,C_1)$-path from $x\in V(C_0)$
to $y\in V(C_1)$ disjoint from $o$. Let $X$ be an $(o,x)$-path on $C_0$, and $Y$
be an $(o,y)$-path on $C_1$.
Then $|P|=l$, and $|Y|\in \{(k-|X|)+l,|X|+l\}$.
\end{claim}
\begin{pf}
Let $C_2:=X\cup Y\cup P$ and $C_3:=(C_0\cup C_1-\widehat{X}\cup \widehat{Y})\cup P$.
If $C_2$ is odd, $|C_2|+|C_3|=2|P|+k+(k+2l)$,
which implies that $C_3$ is also odd, $|C_2|=|C_3|=k+2l$ and $|P|=l$.
In this case, $|Y|=(k-|X|)+l$. If $C_2$ is even,
then $X'\cup Y\cup P$ is an odd cycle, where $X':=C_0-X$.
Similarly, we have $|P|=l$, and $|Y|=(k-|X'|)+l=|X|+l$.
\end{pf}

\medskip

We write $C_0=ox_1x_2\cdots x_{k-1}o$ and $C_1=oy_1y_2\cdots y_{k+2l-1}o$.
For any $x_i\in C_0-\{o\}$ and $(x_i,C_1)$-path $P_i$ from Claim 1, we denote by $x'_i$ the vertex in $V(P_i\cap C_1)$.
Claim 2 implies that
for any $i$, $|P_i|=l$ and $x'_i\in \{y_{l+k-i},y_{l+i}\}$.
\begin{claim}
$l\geq 2$.
\end{claim}
\begin{pf}
Suppose that $l=1$. Set $i:=\frac{k-1}{2}$.
Note that $|P_i|=|P_{i+1}|=1$ (so $P_i,P_{i+1}$ are disjoint) and by Claim 2, $x'_{i}, x'_{i+1}\in \{y_{i+1},y_{i+2}\}$. If $x'_{i}=
x'_{i+1}$, then there is a triangle, a contradiction.
By symmetry, assume that $x'_{i}=y_{i+1}$ and $x'_{i+1}=y_{i+2}$. Then
$(C_1-y_{i+1}y_{i+2})\cup P_i\cup x_{i}x_{i+1}\cup P_{i+1}$ is an odd
cycle of length $k+2l+2$, a contradiction.
\end{pf}

\begin{claim}
$P_1$ and $P_{k-1}$ are disjoint.
\end{claim}
\begin{pf}
Suppose that $P_1$ and $P_{k-1}$ intersect.
Let $z\in V(P_1\cap P_{k-1})$ such that
$|x_{k-1}P_{k-1}z|$ is minimal. Let $R:=x_1P_1z$,
$S:=x_{k-1}P_{k-1}z$ and $T:=zP_1x'_1$.
By Claim 2, $|R|+|T|=l=|S|+|T|$, thus $|R|=|S|$.
Then $C_2:=(C_0-\{o\})\cup R\cup S$ is an odd cycle of length $k-2+2|R|\in \{k,k+2l\}$, implying that $|R|=1$ or $l+1$.
Since $|R|\leq |P_1|=l$, we have $|R|=1$.
If $z\neq x'_1$, then $C_2$ and $C_1$ are disjoint,
a contradiction to Lemma \ref{Lem:2oddcycle-1}.
So $z=x'_1$. Then $R$ is a $(C_0,C_1)$-path, which is disjoint from $o$.
By Claims 2 and 3, $|R|=l\ge 2$, again a contradiction.
\end{pf}

\medskip

By Claims 2 and 4, $\{x'_1,x'_{k-1}\}=\{y_{l+1},y_{k+l-1}\}$.
Let $Y$ be an $(x'_1,x'_{k-1})$-path on $C_1$ with length $2l+2$.
Then $(C_0-\{o\})\cup P_1\cup P_{k-1}\cup Y$ is an
odd cycle of length $k+4l$, a contradiction.
This proves Lemma \ref{Lem:(k+2l)+(k)2}. {\hfill$\Box$}

\subsection{Proof of Lemma \ref{Lem:Two(k)2}.}
We prove by contradiction. Suppose that every two odd cycles
intersect in at most one vertex. By Lemma \ref{Lem:2oddcycle-1}, we
know every two odd cycles intersect in at least one vertex. Thus,
there exist two $k$-cycles $C_0,C_1$ in $G$ with $|V(C_0\cap C_1)|=1$,
and denote the vertex in $C_0\cap C_1$ as $o$.
Write $C_0=ox_1x_2\ldots x_{k-1}o$ and $C_1=oy_1y_2\ldots y_{k-1}o$.

\setcounter{claim}{0}

\begin{claim}
Let $P$ be any $(C_0,C_1)$-path from $x\in V(C_0)$
to $y\in V(C_1)$ disjoint from $o$.
Then $|P|=l$ or $2l$.
\end{claim}

\begin{pf}
We choose $X$ as an $(o,x)$-path on $C_0$, and $Y$ as an $(o,y)$-path
on $C_1$, such that $C_2:=X\cup Y\cup P$ is odd. Thus $C_3:=(C_0\cup C_1-\widehat{X}\cup \widehat{Y})\cup P$ is odd.
Since $|C_2|+|C_3|=2k+2|P|\in \{2k+2l,2k+4l\}$,
we have $|P|\in \{l,2l\}$.
\end{pf}

\begin{claim}
For any $x_i\in C_0-\{o\}$, there is a $P_i$ from $x_i$ to a vertex in $C_1-\{o\}$ (say $x_i'$), internally
disjoint from $C_0\cup C_1$. Similarly, for any $y_j\in C_1-\{o\}$, there is a path $Q_j$
from $y_j$ to a vertex in $C_0-\{o\}$ (say $y_j'$), internally disjoint from $C_0\cup C_1$.
\end{claim}

\begin{pf}
By symmetry, consider an arbitrary vertex $x_i$ in $C_0-\{o\}$. Suppose
that there exists some $C_0$-bridge $H$ such that: $x_i\in V(H-\{o\})$
and $V(H\cap C_1)\subseteq \{o\}$. Let $G_0=:G[H\cup C_0]$ and $G_1:=G-(H-C_0)$ such that $G=G_0\cup G_1$.
Note that $G_1$ is a proper subgraph of $G$ and thus has a proper 3-coloring $\varphi$. By Lemma \ref{Lem:(k+2l)+(k)2}, we know $L(G_0)=\{k\}$.
Then Theorem \ref{THkrsextending} ensures that the restriction $\varphi$ on $C_0$ can be extended to a proper 3-coloring of $G_0$.
Thus $G$ is 3-colorable, a contradiction.
\end{pf}

\medskip

In the following of this subsection, for any vertex $x_i\in C_0-\{o\}$
and any $(x_i,C_1)$-path $P_i$ from Claim 2,
we denote by $x'_i$ the end vertex of $P_i$ in $V(C_1)$. And we also define $y'_j$ for $y_j\in C_1-\{o\}$ analogously.
The next claim summarizes the possible locations of $x'_i$ and $y'_j$, which can be obtained along the same line as in the proof of Claim 2 in Section 5.1.

\begin{claim}
If $|P_i|=2l$, then $x'_i\in \{y_{k-i},y_i\}$; if $|P_i|=l$, then
$x'_i\in \{y_{i-l},y_{i+l},y_{k-i-l},y_{k-i+l}\}$.
Similarly, if $|Q_j|=2l$, then $y'_j\in \{x_j,x_{k-j}\}$; if $|Q_j|=l$,
then $y'_j\in \{x_{j-l},x_{j+l},x_{k-j-l},x_{k-j+l}\}$.

In particular, for each $i\in \{1,k-1\}$, we have the following:
if $|P_i|=2l$, then $x'_i\in \{y_1,y_{k-1}\}$; if $|P_i|=l$,
then $x'_i\in \{y_{1+l},y_{k-l-1}\}$. Similarly,
if $|Q_i|=2l$, then $y'_i\in \{x_1,x_{k-1}\}$; if $|Q_i|=l$,
then $y'_i\in \{x_{1+l},x_{k-l-1}\}$.
\end{claim}

For convenience, we draw $C_0,C_1$ on the plane such that $o,x_1,x_2,...,x_{k-1}$
appear in $C_0$ in the clockwise order, and $o,y_1,y_2,...,y_{k-1}$ appear in
$C_1$ in the counterclockwise order.

\begin{claim}
Let $P_i,P_j$ be $(C_0,C_1)$-paths, from $x_i,x_j\in V(C_0)$ to
$x'_i,x'_j\in V(C_1)$, respectively, where $i<j$.
Let $X$ be the $(x_i,x_j)$-path on $C_0$ not containing $o$.
Then the following hold:\\
(1) If $|P_i|=l$ or $|P_j|=l$, then $P_i,P_j$ are internally disjoint.\\
(2) Suppose that $|P_i|=|P_j|=2l$. If $|X|$ is odd, then $P_i$ and $P_j$ are internally disjoint;
if $\lvert X \rvert$ is odd and $x_i'=x_j'$, then $i=2l$.
In particular, when $\{i,j\}=\{1,2\}$ or $\{1,k-1\}$, $P_i,P_j$ are disjoint.
\end{claim}
\begin{pf}
(1) By symmetry, suppose that $|P_i|=l$ and $P_i,P_j$ intersect on some vertex not in $C_1$.
Let $w\in V(P_i\cap P_j)-V(C_1)$ such that $|wP_jx_j|$ is minimal. Let $P:=x_iP_iw\cup wP_jx_j$,
$C_2:=X\cup P$, and $C_3:=(C_0-\widehat{X})\cup P$. Since $C_2$ and $C_1$ are
disjoint, $C_2$ is even. So $C_3$ is odd, and since
$V(C_1\cap C_3)=\{o\}$, we infer that $|C_3|=k$ by Lemma \ref{Lem:(k+2l)+(k)2}. But $wP_ix'_i$ is a $(C_1,C_3)$-path,
disjoint from $o$, with the length less than $l$, a contradiction to Claim 1.

(2) Suppose that there exists $w\in V(P_i\cap P_j)-V(C_1)$. Choose $w$ such that $|x_jP_jw|$ is minimal.
Let $P=x_iP_iw\cup wP_jx_j$.
If $|P|$ is even, then $X\cup P$ is an odd cycle disjoint from $C_1$, a contradiction to Lemma \ref{Lem:2oddcycle-1}.
So $|P|$ is odd, then $C_2:=P\cup (C_0-\widehat{X})$ is also odd. As $V(C_2\cap C_1)=\{o\}$, we infer that
$|C_2|=k$ by Lemma 12. Note that $wP_{i}x'_i$ is a $(C_2,C_1)$-path with length less than $2l$. Thus $|wP_{i}x'_i|=l$, and
$|x_iP_iw|=|x_jP_jw|=l$. This implies $P$ is even, a contradiction.

Suppose $V(P_i\cap P_j)=\{x'_i\}$. Then $C_3:=X\cup P_i\cup P_j$ is an odd cycle.
As $V(C_3\cap C_1)=\{x'_i\}$, we infer that $|C_3|=|X|+4l=k$.
Note that $oC_0x_i,x_jC_0o$ are two $(C_3,C_1)$-paths disjoint from $x'_i$.
By Claim 1 and Claim 4(1), since $|X|=k-4l$, we deduce $i=|oC_{0}x_i|=|x_jC_0o|=2l$ .

Now let $\{i,j\}=\{1,2\}$. In this case, $|X|=1$ and $|X|$ is odd, so $P_i,P_j$ are internally disjoint.
Suppose that $P_i, P_j$ are not disjoint. Then $V(P_i\cap P_j\cap C_1)=\{x'_i\}=\{x'_j\}$, implying that $i=2l$, a contradiction to $i=1$ (as $i<j$).
The remaining case $\{i,j\}=\{1,k-1\}$ can be proved similarly (as $|X|=k-2$ is also odd).
This proves (2).
\end{pf}

\begin{claim}
$|P_1|=|P_{k-1}|=|Q_1|=|Q_{k-1}|=l$.
\end{claim}
\begin{pf}
Suppose not. By symmetry, assume that $P_1$ is of length $2l$ from $x_1$ to $y_1$,
so we may further assume $x'_1=y_1$ by symmetry and Claim 3.
Note that $P_1$ can also be viewed as $Q_1$.
Suppose that $|P_{k-1}|=2l$ or $|Q_{k-1}|=2l$ (let us say $|P_{k-1}|=2l$).
By Claim 4, $P_1$ is disjoint from $P_{k-1}$, and thus $x'_{k-1}=y_{k-1}$ (because $x'_{k-1}\in \{y_1,y_{k-1}\}$).
Then $(C_1-\{o\})\cup P_1\cup P_{k-1}\cup (x_1ox_{k-1})$
is an odd cycle of length $k+4l$, a contradiction.
So $|P_{k-1}|=|Q_{k-1}|=l$, where
$x'_{k-1}\in \{y_{l+1},y_{k-l-1}\}$ and
$y'_{k-1}\in \{x_{l+1},x_{k-l-1}\}$.

Suppose that $x'_{k-1}=y_{l+1}$ or $y'_{k-1}=x_{l+1}$.
By symmetry, let $P_{k-1}$ be from $x_{k-1}$ to $y_{l+1}$. Then
$P_1\cup x_1C_0x_{k-1}\cup y_{l+1}C_1y_1\cup P_{k-1}$
is an odd cycle of length $k+4l-2$, implying $l=1$.
So $P_{k-1}$ is from $x_{k-1}$ to $y_2$.
If $Q_{k-1}=y_{k-1}x_2$, then $(x_1oy_{k-1})\cup Q_{k-1}\cup x_{2}C_0x_{k-1}\cup P_{k-1}\cup y_2y_1\cup P_1$
is an odd cycle of length $k+2l+2$, a contradiction.
So $Q_{k-1}=y_{k-1}x_{k-2}$, but then
$x_1C_0x_{k-2}\cup Q_{k-1}\cup (y_{k-1}ox_{k-1})\cup P_{k-1}\cup y_2y_1\cup P_1$
is an odd cycle of length $k+2l+2$, again a contradiction.
Hence we may assume that $P_{k-1}$ is from $x_{k-1}$ to $y_{k-l-1}$,
and $Q_{k-1}$ is from $y_{k-1}$ to $x_{k-l-1}$.

If $k\neq l+2$, then $y_{k-l-1}\neq y_1$.
Let $X$ be the $(y_1,y_{k-l-1})$-path on $C_1$ containing $o$ and with $|X|=l+2$.
Then $P_1\cup X\cup P_{k-1}\cup (C_0-\{o\})$ is
an odd cycle of length $k+4l$, a contradiction.
Thus $k=l+2$, and now $P_{k-1}$ is from $x_{k-1}$ to $y_1$,
and $Q_{k-1}$ is from $y_{k-1}$ to $x_1$.
Recall that $l=k-2$, implying that $l\ge 3$ is odd.
Consider the path $P_2$ from Claim 2.
If $|P_2|=l$, then by Claim 3, $x'_2\in \{y_{2-l},y_{l+2},y_{k-2-l},y_{k-2+l}\}$, contradicting the facts that $k=l+2$ and $l\ge 3$.
So $|P_2|=2l$. Then $x'_2\in \{y_2,y_{k-2}\}$, and $P_2$ is internally disjoint with $P_1$ or $P_{k-1}$ (by Claim 4).
If $x'_2=y_{k-2}$, then $P_2\cup (C_1-\{o,x_1\})\cup P_{k-1}\cup (C_0-\{o,y_{k-1}\})$ is an odd cycle of length $k+4l-4>k+2l$ (as $l\ge 3$), a contradiction.
So $x'_2=y_2$. Then $(C_1-y_1y_2)\cup P_1\cup x_1x_2\cup P_2$
is an odd cycle of length $k+4l$, a contradiction.
The proof of this claim is complete.
\end{pf}

\medskip

By Claims 4 and 5, we see that for any distinct $i,j$, where $i\in \{1,k-1\}$,
\begin{equation}\label{equ:ij}
P_i, P_j \text{ (and respectively, } Q_i, Q_j) \text{ are internally disjoint.}
\end{equation}

\begin{claim}
For any $i,j\in \{1,k-1\}$, $P_i$ and $Q_j$ are disjoint.
\end{claim}
\begin{pf}
By symmetry, it will suffice to show that $P_1$ and $Q_1$ are disjoint. Suppose for a contradiction that $P_1,Q_1$ are not disjoint.

We first show that $P_1$ can be chosen from $x_1$ to $y_1$ and $k=l+2$ (thus $l\ge 3$ is odd).
Since $|Q_1|=l$, by Claim 3, we have $y_1'\in \{x_{1+l},x_{k-l-1}\}$.
If $y_1'=x_1$, then we must have $k=l+2$ and one can view the $(y_1,x_1)$-path $Q_1$ as $P_1$, done.
So $y_1'\neq x_1$ and thus $Q_1$ can be viewed as some $P_j$ for $j\neq 1$.
By \eqref{equ:ij}, $P_1,Q_1$ (which are viewed as $P_1,P_j$) are internally disjoint, but not disjoint.
So we have either $x'_1=y_1$ or $y_1'=x_1$.
By symmetry (as we shall see) assume that $x'_1=y_1$. Since $|P_1|=l$, by Claim 3, $x'_1\in \{y_{1+l},y_{k-l-1}\}$, which forces $x'_1=y_1=y_{k-l-1}$. So $k=l+2$ and $P_1$ is from $x_1$ to $y_1$.

Next we claim that $P_1, P_{k-1}$ are disjoint, where $P_{k-1}$ is from $x_{k-1}$ to $y_{k-1}$. By \eqref{equ:ij}, $P_1,P_{k-1}$ are internally disjoint and by Claims 5 and 3, $x'_{k-1}\in \{y_1,y_{k-1}\}$.
If $x'_{k-1}=y_1$, then $(C_0-\{x_{k-1}o, ox_1\})\cup P_1\cup P_{k-1}$ is an odd cycle of length $k+2l-2$, which implies that $l=1$ and $k=3$, a contradiction to $k\ge 5$.
So $x'_{k-1}=y_{k-1}$ and thus $P_1,P_{k-1}$ are disjoint.

In the following, we will consider $P_2, P_{k-2}$ to find an odd cycle of length larger than $k+2l$, which is a contradiction.
If $|P_2|=l$, in view of the facts $k=l+2$ and $l\ge 3$, there is no valid choice for $x_2'\in \{y_{2-l},y_{2+l},y_{k-2-l},y_{k-2+l}\}$ according to Claim 3.
So $|P_2|=2l$ and by Claim 3, $x'_{2}\in \{y_2,y_{k-2}\}$.
Similarly, we have $|P_{k-2}|=2l$ and $x'_{k-2}\in \{y_2,y_{k-2}\}$.
Suppose that $x'_2=x'_{k-2}$ (which is either $y_2$ or $y_{k-2}$).
Let $X$ be the $(x_2,x_{k-2})$-path on $C_0$ avoiding $o$. So $|X|=k-4$ is odd.
Then $X\cup P_2\cup P_{k-2}$ is an odd cycle of length $k+4l-4>k+2l$ (as $l\ge 3$), a contradiction.
So we have $x'_2\neq x'_{k-2}$ and thus $P_1,P_2,P_{k-2},P_{k-1}$ are pairwise disjoint.
In either case of $x'_2\in \{y_2,y_{k-2}\}$, we see that $P_1\cup y_1y_2\cup P_2\cup (C_0-\{x_1x_2,x_{k-1}x_{k-2}\})\cup P_{k-2}\cup y_{k-2}y_{k-1}\cup P_{k-1}$
is an odd cycle of length $k+6l>k+2l$, a contradiction. This proves Claim 6.
\end{pf}

\begin{claim}
$P_1,P_{k-1}$ share the endpoint in $C_1$, or $Q_1,Q_{k-1}$ share the endpoint in $C_0$.
\end{claim}
\begin{pf}
Otherwise, $P_1,P_{k-1},Q_1,Q_{k-1}$ are pairwise disjoint.
By symmetry, assume that $y_1'\in y_{k-1}'C_0o$.
Let $X$ be the $(x_1',x_{k-1}')$-path on $C_1$ not containing $o$.
Then $x_1C_0y_{k-1}'\cup Q_{k-1}\cup (y_{k-1}oy_1)\cup Q_1\cup y_1'C_0x_{k-1}\cup P_{k-1}\cup X\cup P_1$
is an odd cycle of length $k+4l$, a contradiction.
\end{pf}

\begin{claim}
$l=1$, $N(x_1)\cap N(x_{k-1})\cap \{y_2,y_{k-2}\}\neq \emptyset$ and $N(y_1)\cap N(y_{k-1})\cap \{x_2,x_{k-2}\}\neq \emptyset.$
\end{claim}
\begin{pf}
By Claim 7, assume by symmetry that $Q_1, Q_{k-1}$ are from $y_1,y_{k-1}$ to $x_{l+1}$ respectively.
Then $C_2:=Q_1\cup Q_{k-1}\cup (C_1-\{o\})$ is an odd cycle
of length $k-2+2l$, which intersects $C_0$ only on $x_{l+1}$.
By Lemma \ref{Lem:(k+2l)+(k)2}, $|C_2|=k$ and thus $l=1$.
Suppose $P_1,P_{k-1}$ can be chosen to be disjoint, say $P_1=x_1y_{k-2}$ and $P_{k-1}=x_{k-1}y_2$. But then the cycle
$(C_0-\{o\})\cup P_1\cup (y_{k-2}y_{k-1}oy_1y_2)\cup P_{k-1}$ is an odd cycle of $k+2l+2$,
a contradiction. This proves the claim.
\end{pf}

\medskip

By symmetry, we may assume that $P_1=x_1y_2$, $P_{k-1}=x_{k-1}y_2$,
$Q_1=y_1x_2$ and $Q_{k-1}=y_{k-1}x_2$. Let $C_2:=(C_0-\{o\})\cup (x_1y_2x_{k-1})$.
Note that $|C_2|=|C_1|=k$ and $V(C_2\cap C_1)=\{y_{2}\}$.
We then can treat $C_2,y_2$ as the new $C_0,o$, and thus  all previous claims hold for $C_1$ and $C_2$.
In particular, by Claim 8, we have $N(y_1)\cap N(y_3)\cap \{x_2,x_{k-2}\}\neq \emptyset$.
If $y_1x_2,y_3x_2\in E(G)$, as $Q_{k-1}=y_{k-1}x_2$,
$G$ has an odd cycle $(C_1-\{o,y_1,y_2\})\cup (y_{k-1}x_2y_3)$ of length $k-2$;
otherwise $y_1x_{k-2},y_3x_{k-2}\in E(G)$, then, as $Q_1=y_1x_2$,
$G$ has an odd cycle $(C_0-\{o,x_1,x_{k-1}\})\cup (x_{k-2}y_1x_2)$ of length $k-2$, a contradiction.
Lemma \ref{Lem:Two(k)2} now is proved.
This, together with Lemmas \ref{Lem:Two(k+2l)} and \ref{Lem:(k+2l)+(k)2}, complete the proof of Lemma \ref{Lem:2oddcycle}.
{\hfill$\Box$}

\medskip

\bigskip

\noindent {\bf Acknowledgements.} The authors would like to thank Rong Luo for pointing out the question,
and an anonymous referee for carefully reading the manuscript and for giving valuable
comments which helped improving the presentation of the paper.

\end{document}